\PassOptionsToPackage{dvipsnames}{xcolor}
\documentclass[leqno]{siamart1116}
\usepackage{adjustbox}
\usepackage{amssymb}
\usepackage{bm,bbm}
\usepackage{caption}
\usepackage{color,colortbl}
\usepackage{dsfont}
\usepackage{dutchcal}
\usepackage{enumerate}
\usepackage{enumitem}
\usepackage{graphicx}
\usepackage{mydef}
\usepackage{marginnote}
\usepackage{multirow}
\usepackage[normalem]{ulem}
\usepackage[numbers]{natbib}
\usepackage{subcaption}
\usepackage{xspace}
\usepackage{dashrule}

\DeclareRobustCommand\full{\hdashrule[0.5ex][c]{6.4mm}{0.6pt}{}}
\DeclareRobustCommand\dashed{\hdashrule[0.5ex][c]{6.8mm}{0.6pt}{1.9mm 0.6mm 0.3mm 0.6mm}\!}
\DeclareRobustCommand\dotted{\hdashrule[0.5ex][x]{6.4mm}{0.6pt}{0.9mm}\!}
\DeclareRobustCommand\fulldot{\hdashrule[0.5ex][x]{6.4mm}{0.6pt}{0.3mm}\!}

\let\cite=\citet
\setlist{leftmargin=20pt}

\begin{document}

\newcommand{\TheTitle}{Preconditioning of the generalized Stokes problem
  arising from the approximation of the time-dependent Navier-Stokes equations.}
\newcommand{\TheAuthors}{Melvin Creff, J.-L. Guermond}

\renewcommand{\thefootnote}{\arabic{footnote}}

\headers{Solving the Schur complement}{\TheAuthors} 

\title{{\TheTitle}%
  \thanks{Draft version, \today%
    \funding{This material is based upon work supported in part by the
      National Science Foundation, under grant DMS-2110868 (JLG), by
      the Air Force Office of Scientific Research, USAF, under
      grant/contract number FA9550-23-1-0007 (JLG), the Army Research
      Office, under grant number W911NF-19-1-0431 (JLG), and the
      U.S. Department of Energy by Lawrence Livermore National
      Laboratory under Contracts B640889, B641173 (JLG).}}}

\author{Melvin Creff\footnotemark[2]
  \and 
  Jean-Luc~Guermond\footnotemark[2]}

\maketitle

\renewcommand{\thefootnote}{\fnsymbol{footnote}}

\footnotetext[2]{%
  Department of Mathematics, Texas A\&M University, 3368 TAMU,
  College Station, TX 77843, USA.}

\renewcommand{\thefootnote}{\arabic{footnote}}

\begin{abstract}
  The paper \bal{compares} standard iterative methods for solving the
  generalized Stokes problem arising from the time and space
  approximation of the time-dependent incompressible Navier-Stokes
  equations. Various preconditioning techniques are considered (Schur
  complement, fully coupled system, with and without augmented
  Lagrangian). One investigates whether these methods can compete with
  traditional pressure-correction and velocity-correction methods in
  terms \bal{of throughput (number of degrees of freedom per time step
    per core per second).}  Numerical tests on fine unstructured meshes (68
  millions degrees of freedoms) demonstrate GMRES/CG convergence rates
  that are independent of the mesh size and improve with the Reynolds
  number \bal{for most methods.  Three conclusions are drawn: (1)
    Whether solving the pressure Schur complement or the fully coupled
    system does not make any significant difference in terms of
    throughput. (2) Although very good parallel scalability is
    observed for the augmented Lagrangian method, the best throughput
    is achieved without using the augmented Lagrangian
    formulation. (3) The throughput of all the methods tested in the
    paper are on average 25 times slower than that of traditional
    pressure-correction and velocity-correction methods. Hence,
    although all these methods are very efficient for solving steady
    state problems, none of them is unfortunately competitive for
    solving time-dependent problems.}
\end{abstract}



\pagestyle{myheadings}
\thispagestyle{plain}
\markboth%
  {M. Creff, J.-L. Guermond}%
  {Augmented Lagrangian}


\section{Introduction}
\bal{By far,} the most popular time-stepping techniques for solving the
time-dependent Navier-Stokes equations are the so-called projection
methods initially proposed by Chorin~\citep{Chor68} and
Temam~\citep{tema69}. Although significant progresses have been made
since the original work of Chorin and Temam (see \eg
\citep{Guermond_Minev_Shen_2006} for a brief overview of the
literature on the topic), all these methods have however limited
accuracy in time. Moreover, when combined with space discretization,
they do not allow for a good control of the weak divergence of the
velocity field. A better control on the accuracy and on the weak
divergence of the velocity can be obtained by using more traditional
time-stepping techniques combining implicit and explicit Runge-Kutta
techniques. The downside of these techniques though is that one has to
solve at each time step a generalized Stokes problem with a saddle
point structure. Many iterative techniques to solve this problem have
been proposed in the literature (see \eg
\cite[Chap.~8]{Elman_Silvester_book_2005} for a review of the
literature). They all more or less \bal{involve the pressure Schur
  complement at some point. But} the
condition number of the Schur complement grows unboundedly as the time
step and the viscosity go to zero. This in turn implies that the
performance of these methods strongly depends on how they are
preconditioned.  \bal{The objective of the paper is to revisit
  preconditioning methods for solving the generalized
  Stokes problem arising from the time approximation of the Navier-Stokes equations, and
  determine whether they can match projection methods in terms of
  throughput.}

The solution of the generalized Stokes problem has a long history that
started with the landmark papers of \cite{Cahouet_Chabard_1988} and
\cite{Fortin_Glowinski_Book_1983}.  In \citep{Cahouet_Chabard_1988},
the authors propose a preconditioner of the pressure Schur complement
that behaves well in the two limits $\mu\dt h^{-2} \to 0$ and
$\mu\dt h^{-2} \to \infty$ where $\mu$ is the viscosity, $\dt$ is the
time step, and $h$ is the mesh size. In
\citep{Fortin_Glowinski_Book_1983} the authors introduce the augmented
Lagrangian method and show that the conditioning of the pressure Schur
complement improves as the scalar multiplier of the augmented
Lagrangian term increases (\cf \eg
\cite[Prop.~2.1]{Golub_Greif_SISC_2003}
\cite[Lem.~4.1]{Benzi_Olshanksii_SISC_2006} and
Proposition~\ref{prop:schur_augmented_lagrange} below).  \bal{An
  alternative approach to the Uzawa iterations on the pressure Schur
  complement suggested in
  \cite{Rusten_Winther_SIAM_Mat_Anal_Appl_1992},
  \cite{Elman_Silvester_SISC_1996}, consists of iteratively solving
  the coupled velocity-pressure system. This method also involves the
  pressure Schur complement, but its main advantage is that it does
  not require to solve the velocity problem with very high precision
  when estimating the action of the Schur complement.}  This method is
used with the augmented Lagrangian in
\cite{Benzi_Olshanksii_SISC_2006},
\cite{Benzi_Olshanskii_Wang_IJNMF_2011}, and
\cite{Benzi_Wang_SISC_2011}. In these papers the authors show that
combining the Cahouet\&Chabard preconditioner and the augmented
Lagrangian technique gives a good preconditioning method of the
generalized Stokes problem.  This method has been numerically
evaluated in many papers like
\cite{Farrell_Mitchell_Wechsung_SISC_2019},
\cite{Moulin_Jolivet_Marquet_CMAME_2019}, and
\cite{Shih_Stadler_Wechsung_SISC_2023}, where it is indeed shown that
the number of MINRES or GMRES iterations used to solve the pressure
Schur complement problem behaves well with respect to the mesh size
and the viscosity.

\bal{We revisit the preconditioning methods discussed above (for the
  pressure Schur complement and the fully coupled system) assuming
  that the generalized Stokes problem is generated from the time and
  space approximation of the Navier-Stokes equations with the
  nonlinear term made explicit. The time step is set to satisfy the
  standard CFL restriction, see \eqref{def_of_dt}. The approximation
  in space is done with mixed finite elements.}
Thorough numerical tests on large two-dimensional problems are done
\bal{on a HPC cluster using PETSc for the parallel linear algebra} and
mixed Hood--Taylor finite elements $\polP_2/\polP_1$ and
$\polP_3/\polP_2$. The tests are also done with the generic Newtonian
form of the viscous tensor inducing nontrivial couplings between the
Cartesian components of the velocity. This allow us to test the
robustness of the methods and their capability to handle nonlinear
expressions for the stress tensor.

Three key conclusions are drawn from the paper. The first one is that
whether solving the pressure Schur complement problem or the fully
coupled system does not make any significant difference in term of
throughput \bal{(see Figure~\ref{fig:efficiency_best}).} The second
conclusion is that although very good parallel scalability is observed
for the augmented Lagrangian formulation,
the best throughput is achieved without the augmented Lagrangian
term. The two main reasons are that \bal{it takes too many inner
  iterations to solve the velocity problem associated with the
  augmented Lagrangian term when using the pressure Schur complement
  method, and it takes too many outer GMRES iterations when
  using the fully coupled method.
  Note that this conclusion only holds true in the context the
  generalized Stokes problem arising from the time approximation of
  the Navier-Stokes equation.  The augmented Lagrangian formulation is
  very efficient for solving steady state problems.
  The third conclusion, summarized in Table~\ref{tab:litterature_comparisons},}
is that the throughput of all the methods tested in the
paper is nowhere near that of projection methods. They are on average
25 times slower than traditional pressure-correction and
velocity-correction methods.
Hence, although all these methods are
very efficient and close to optimality for solving steady state
problems, our conclusion, at the time of this writing, is that they
cannot reasonably compete with projection methods or variations
thereof for solving the time-dependent Navier--Stokes equations.
\bal{This paper demonstrate that to achieve a throughput comparable to
  projection methods while being more robust and accurate, new
  solution methods accounting for both time and space must be
  invented.  This paper is an invitation to start new
  lines of investigations in this direction.}

The paper is organized as follows. First we formulate the problem in
\S\ref{Sec:Formulation}. We describe the discrete setting and review
preconditioners combining ideas from \cite{Cahouet_Chabard_1988} and
\cite{Fortin_Glowinski_Book_1983} further developed in
\cite{Benzi_Olshanksii_SISC_2006}, \citep{Benzi_Wang_SISC_2011},
\cite{Benzi_Olshanskii_Wang_IJNMF_2011},
\cite{Elman_Silvester_SISC_1996},
\cite{Rusten_Winther_SIAM_Mat_Anal_Appl_1992}.  We discuss in
\S\ref{Sec:Implementation} and
\S\ref{Sec:solving_pressure_schur_complement} key practical
difficulties that arise when implementing these preconditioners. We
particularly focus on the parallel aspects of the question.  In
\S\ref{Sec:Implementation} we discuss and tests preconditioners to
solve the velocity problem. We focus on the couplings induced either by
the augmented Lagrangian term or the Newtonian form of the viscous
tensor.
\bal{Discarding}
every coupling block in the velocity matrix is the best
preconditioning strategy (see preconditioner $\tA_{\lambda,3}$ defined
in \eqref{precond_A_plus_lambda_D_3}). In
\S\ref{Sec:solving_pressure_schur_complement} we discuss
preconditioning techniques for the pressure Schur complement.  All
these preconditioners involve solving $B M_\bV^{-1} B\tr \sfX =\sfY$,
where $B$ is the matrix associated with the divergence operator and
$M_\bV$ is the velocity mass matrix. Various preconditioning for
$B M_\bV^{-1} B\tr$ are tested in
\S\ref{Sec:solving_pressure_schur_complement}.  In
\S\ref{Sec:Numerical_illustrations}, we test preconditioners for the
pressure Schur complement. Robustness with respect to the viscosity,
the time step, the mesh size, and the value of the coefficient scaling
the augmented Lagrangian is thoroughly investigated.  In
\S\ref{Sec:Numerical_performance_Method2} we test preconditioning
techniques for the full generalized Stokes system.  \bal{We
  summarize our conclusions in \S\ref{Sec:Conclusions}. (1) Whether
  iteratively solving the pressure Schur complement or the fully
  coupled system does not make any significant difference in terms of
  throughput (see Figure~\ref{fig:efficiency_best}). (2) The best
  throughput is achieved without using the augmented Lagrangian
  formulation. (3) The throughput of all the methods tested in the
  paper are on average 25 times slower than that of traditional
  pressure-correction and velocity-correction methods (see
  Tables~\ref{tab:litterature_comparisons}-\ref{tab:pred_corr_eff}).}

\section{Formulation of the problem} \label{Sec:Formulation}
We formulate the problem and introduce the notation in this section.

\subsection{Semi-discrete problem}
We consider the generalized Stokes problem arising from the
approximation in time of the time-dependent Navier-Stokes equations in
a domain $\Dom\subset \Real^d$.  Here $\Dom$ is an open bounded
Lipschitz polyhedron. We denote by $\ell_\Dom$ the diameter of $\Dom$.
The exact form of the time discretization that is adopted to
approximate the problem in time \bal{is not the concern of this paper and could be
any implicit method like (BDF, IRK, etc.)}. Let
$\bV$ be a closed subspace of $\bH^1(\Dom)$ and $Q$ be a closed
subspace of $L^2(\Dom)$. All these spaces are equipped with their natural
norm $\|\cdot\|_{\bH^1(\Dom)}$, $\|\cdot\|_{L^2(\Dom)}$. We assume
that we are given at each time step a linear form
$\bef\in \calL(\bV;\Real)$, and the problem is reduced to seeking a
pair $(\bu,p)\in\bV\CROSS Q$ so that
\begin{equation}
  a(\bu,p) - b(\bv,p) + b(\bu,q) = \bef(\bv),\qquad
  \forall (\bv,q)\in\bV\CROSS Q, \label{weak_Stokes} 
\end{equation}
where the bilinear form $a$ and $b$ are defined by
\begin{equation}
 a(\bv,\bw) \eqq \int_\Dom \left(\frac{1}{\dt}\bv\SCAL \bw +
 2\mu\pole(\bv){:}\pole(\bw)\right) \diff x,\quad
 b(\bv,q) \eqq \int_\Dom q \DIV \bv\diff x.
\end{equation}
Here $\pole(\bv)\eqq\frac12(\GRAD \bv+(\GRAD\bv)\tr) $ is the strain
rate tensor, $\mu$ is the shear viscosity, and $\dt$ is the time
step \bal{normalized to includes other factors from the time-steeping method}.
The boundary conditions of the problem are encoded in the
definition of the space $\bV$ and the linear form $\bef$.  This
problem has a unique solution (see \eg \cite{GR86},
\cite[\S~4.2]{Boffi_Brezzi_Fortin_book_2013}).  Our objective is to
construct an iterative method for the solution of a discrete version
of \eqref{weak_Stokes} that performs well with respect to the product
$\dt\mu$ and the mesh size.

\subsection{Discrete problem}

In the rest of this paper we solely focus on a discrete version of
\eqref{weak_Stokes}. We consider two sequences of finite-dimensional
vector spaces $\{\bV_h\}_{h\in\calH}$ and $\{Q_h\}_{h\in\calH}$ where
the index set $\calH$ is countable and has $0$ as unique accumulation
point.  The velocity is approximated in $\bV_h$ and the pressure is
approximated in $Q_h$. To simplify some arguments, we assume that the
velocity approximation is conforming, \ie $\bV_h\subset \bV$, and the
pair $(\bV_h,Q_h)$ is inf-sup stable, \ie
\begin{equation}
  \beta_h\eqq \inf_{q_h\in Q_h{\setminus}\{0\}} \sup_{\bv_h\in
    \bV_h{\setminus}\{\bzero\}} \frac{b(\bv_h,q_h)}{\|\bv_h\|_{\bH^1(\Dom)}\|q_h\|_{L^2(\Dom)}}>0.
\end{equation}

The discrete version of \eqref{weak_Stokes} consists of seeking
$\bu_h\in \bV_h$ and $p_h\in Q_h$ so that the following holds true
for all $\bv_h\in \bV_h$ and all $q_h\in Q_h$:
\begin{equation}
 a(\bu_h,\bv_h) - b(\bv_h,p_h) = \bef(\bv_h), \qquad
b(\bu_h,q_h) = 0. \label{discrete_weak_Stokes}
\end{equation}

Let $\{\bphi_i\}_{i\in\calV}$ and $\{\psi_k\}_{k\in\calQ}$ be bases of
$\bV_h$ and $Q_h$ (say, finite-element-based global shape
functions). We define the velocity mass matrix $M_\bV$, the velocity
stiffness matrix $E_\bV$, and the pressure mass matrix $M_Q$ by setting
\begin{subequations}\begin{align}
  M_{\bV\!,ij}&\eqq \int_\Dom \bphi_i\SCAL\bphi_j\diff x,\qquad \forall i,j\in \calV,\\
M_{Q,ij}&\eqq \int_\Dom \bpsi_i\SCAL\bpsi_j\diff x,\qquad \forall i,j\in \calV,\\
E_{\bV\!,kl} &\eqq\int_\Dom2\pole(\bphi_k){:}\pole(\bphi_l) \diff x \qquad \forall k,l\in \calQ.\label{def_of_EV}
\end{align}  \end{subequations}
To formulate the algebraic version of \eqref{discrete_weak_Stokes}, we
introduce the matrices associated with the bilinear forms $a$ and $b$:
\begin{align}
  A_{ij}&\eqq \frac{1}{\dt} M_{\bV\!,ij}+ \mu E_{\bV\!,ij}=a(\bphi_j,\bphi_i),\qquad \forall i,j\in\calV,\\
 B_{kj}&\eqq b(\bphi_j,\psi_k), \qquad \forall k\in\calQ,\ \forall j\in\calV.\label{def_of_Aij_Bij}
\end{align}
\bal{To be representative of situations corresponding to the approximation
of the time-dependent Navier-Stokes equations, where the time step
decreases like the mesh size due to the nonlinearities being made
explicit in time,  in the rest of this paper we set
\begin{equation}
  \tau=N^{-\frac1d},\label{def_of_dt}
\end{equation}
where $N$ is the total number of degrees of freedom for one scalar component
of the velocity and $d$
is the space dimension.  $N^{-\frac1d}$ scales like the mesh size. The
tests reported in the paper are done in two dimensions}.

Denoting $\bu_h\eqq \sum_{i\in\calV}\sfU_i\bphi_i$ and
$p_h\eqq \sum_{k\in\calQ}\sfP_k\psi_k$ the approximations of the
velocity and the pressure, the discrete
problem~\eqref{discrete_weak_Stokes} is then equivalent to solving the
following linear system:
\begin{equation}
  A \sfU -B\tr \sfP = \sfF,\qquad B\sfU= 0. \label{Saddle_point_A_BT_B}
\end{equation}

Two different strategies are advocated in the literature to solve the
saddle point problem~\eqref{Saddle_point_A_BT_B}.  The first one,
which is the oldest, consists of constructing the pressure Schur
complement of the problem and solving it. This method solves first for
the pressure then solves for the velocity.  \bal{We henceforth refer
  to this approach and variations thereof as Method~1.}  The second
approach consists of solving for both the velocity and the pressure at
the same time; see \eg \cite{Rusten_Winther_SIAM_Mat_Anal_Appl_1992},
\cite{Elman_Silvester_SISC_1996},
\cite{Furuichi_May_Tackley_JCP_2021},
\cite{Shih_Stadler_Wechsung_SISC_2023}. \bal{We henceforth refer to
  this approach and variations thereof as Method~2.} We now detail these
two methods.

\subsection{Method~1: The pressure Schur complement}

As the coercivity of $a$ implies that the matrix $A$ is invertible, a
standard way of solving \eqref{Saddle_point_A_BT_B} consists of
constructing the Schur complement with respect to the
pressure. Letting $S\eqq B A^{-1} B\tr $,  the problem \eqref{Saddle_point_A_BT_B}
is equivalent to finding 
$\sfP$ and $\sfU$ so that
\begin{equation}
S \sfP = - B A^{-1}\sfF,\qquad A\sfU = \sfF + B\tr \sfP.
\end{equation}
Note that here the pressure and the velocity problems are somewhat
uncoupled. One first solves for the pressure, then one computes the
velocity.  The pressure problem $S \sfP = - B A^{-1}\sfF$ is in
general solved iteratively.  As the convergence rate of iterative
methods is controlled by the condition number, we now recall what is
known on the $\ell^2$-condition number of $S$.  The following result
can be found in \cite[Prop.~5.24]{Elman_Silvester_book_2005} (see also
\citep[Prop.~50.14]{Ern_guermond_FE_II_2021} for other details on
the topic).
\begin{proposition}[Spectrum of $S$] Let $\sigma(S)$ be the spectrum
  of $S$. Let $\mu_{\min}$, $\mu_{\max}$ be the smallest and largest
  eigenvalues of the pressure mass matrix $M_Q$. Let $\|a\|$ and
  $\|b\|$
  be the norms of the bilinear forms $a$ and $b$. Let $\alpha_h$ be
  the coercivity constant of $a$. Then
  $\sigma(S) \subset [\mu_{\min} \frac{\beta_h^2}{\|a\|} , \mu_{\max}
  \frac{\|b\|^2}{\alpha_h}]$.
\end{proposition}

As $\|a\|\sim \frac{\mu}{\ell_\Dom^2 +\frac{1}{\dt}}$ and
$\alpha_h\sim \frac{\mu}{\ell_\Dom^2}$, the $\ell^2$-condition number
of $S$ scales like
$(1+\frac{\ell_\Dom^2}{\dt \mu})\frac{\|b\|}{\beta_h}$. This number is
larger than $(1+\frac{\ell_\Dom^2}{\dt \mu})$ and grows unboundedly as
$\frac{\ell_\Dom^2}{\dt \mu}\to \infty$. Hence, if not properly
preconditioned, iterative methods for solving
$S \sfP = - B A^{-1}\sfF$ have convergence rates that degrade as the
viscosity decreases and the time step goes to zero.

The loss of convergence rate mentioned above can
be mitigated by using the augmented Lagrangian technique introduced by
\cite{Fortin_Glowinski_Book_1983}.  This method replaces the original
system \eqref{Saddle_point_A_BT_B} by the following equivalent system:
\begin{subequations}\begin{align}
A_{\lambda} &\sfU - B\tr \sfP = \sfF, \qquad B\sfU= 0, \label{LA_problem}\\
A_{\lambda} &\eqq  A + \lambda \mu B\tr M_Q^{-1} B,   \label{def_of_Alambda}                        
\end{align}\end{subequations}
where $\lambda$ is a non-dimensional positive parameter. Here again, one can solve the
problem by constructing the pressure Schur complement,
\begin{subequations}\begin{align}
  S_{\lambda} &\sfP = - B (A + \lambda \mu B\tr M_Q^{-1}B)^{-1}\sfF,\qquad
A\sfU = \sfF + B\tr \sfP,
\label{Schur_LA_problem}\\
  S_{\lambda}&\eqq B (A + \lambda \mu B\tr M_Q^{-1} B)^{-1} B\tr. \label{def_Schur_LA_problem}
\end{align}\end{subequations}
The main properties of this method are summarized in the following
standard result (\cf \eg \cite[Prop.~2.1]{Golub_Greif_SISC_2003} \cite[Lem.~4.1]{Benzi_Olshanksii_SISC_2006})
\begin{proposition}[Augmented Lagrangian]  \label{prop:schur_augmented_lagrange}
  Let $s_\flat$ and $s_\sharp$ be the smallest and largest eigenvalues
 of $S$. Then the following holds true:
  \begin{subequations}\begin{align}
   &S_{\rho}^{-1} = \rho M_Q^{-1} + S^{-1}, \label{eq1:prop:schur_augmented_lagrange}\\
&\sigma(M_Q^{-1} S_{\rho}) \subset [(\rho + s_\flat^{-1})^{-1}, (\rho + s_\sharp^{-1})^{-1}].
\end{align}   \end{subequations}
\end{proposition}

The above estimate of the spectrum of $S_{\lambda}$ shows that the
$\ell^2$-condition number of $S_{\lambda}$ converges to unity as
$\lambda\mu\dt/\ell_\Dom^2 \to \infty$. But, as the condition number
of $S_{\lambda}$ converges to unity, the condition number of
$A_\lambda$ becomes unrealistically large when
$\lambda\mu\dt/\ell_\Dom^2 \to \infty$. Hence, one important
difficulty with the augmented Lagrangian method is to find a value of
$\lambda$ that is balanced.

In the rest of the paper we discuss \bal{left} preconditioning techniques for
$S_{\lambda}$ and the pros and cons for the augmented Lagrangian
method.

\subsection{Method~2: coupled problem}\label{Sec:Method_2}

As observed in \cite{Rusten_Winther_SIAM_Mat_Anal_Appl_1992} (second
to last paragraph of the introduction therein),  for
Krylov methods based on the pressure Schur complement to
work properly, it is necessary to
compute $S\sfR$ exactly for all pressure vector $\sfR$.
This entails that it is imperative to solve the velocity problem
$A\bW = \bK$ exactly for all velocity vector $\sfK$. One way to avoid
this constraint consists of using Krylov methods based on the
full velocity-pressure problem as advocated in
\citep{Rusten_Winther_SIAM_Mat_Anal_Appl_1992}; (see also Eq.~(2.5)
and (2.12) in \cite{Elman_Silvester_SISC_1996}).  The idea is to
exploit the following exact factorization
  \begin{align}
  \polA_\lambda^{-1}\eqq  \begin{pmatrix} A_\lambda & -B\tr\\ B & 0\end{pmatrix}^{-1}
= \begin{pmatrix} I_\bV & A_\lambda^{-1} B\tr\\ 0 & I_Q\end{pmatrix}
\begin{pmatrix} A_\lambda^{-1} & 0 \\ 0 & S_\lambda^{-1}\end{pmatrix}
\begin{pmatrix}I_\bV & 0 \\ -B A_\lambda^{-1} & I_Q\end{pmatrix},\label{exact_factorization}
    \end{align}
    where $I_\bV$ and $I_Q$ are the identity matrices.  This identity
    shows that is suffices to replace $A^{-1}$ and $S^{-1}$ by
    preconditioners to construct preconditioners of $\polA_\lambda$.
    \bal{Denoting by $\tA_\lambda^{-1}$, $\tS_\lambda^{-1}$
      preconditioners of $\tA_\lambda$, $\tS_\lambda$, the second
      method we consider (henceforth called Method~2) consists of
      using GMRES with the following preconditioner:
\begin{align}
\polA_\lambda^{-1}\sim \begin{pmatrix} I_\bV & \tA_\lambda^{-1} B\tr\\ 0 & I_Q\end{pmatrix}
\begin{pmatrix} \tA_\lambda^{-1}  & 0 \\ 0 & \tS_\lambda^{-1}\end{pmatrix}
\begin{pmatrix}I_\bV & 0 \\ -B \tA_\lambda^{-1} & I_Q\end{pmatrix}.\label{inexact_factorization_of_Alambda}
\end{align}}%
No exact inverse has to be computed. It should also be noted that,
contrary to appearances, this method only requires two inversions of
$A_\lambda$ at each Krylov iteration.

\bal{In the rest of the paper we are going to investigate left
  preconditioning techniques for $\polA_\lambda$ using variations of the preconditioners $\tA_\lambda^{-1}$ and $\tS_\lambda^{-1}$ 
  in \eqref{inexact_factorization_of_Alambda}.  We are also going to investigate
  whether the augmented Lagrangian version of the method is indeed
  superior to the non augmented one ($A_\lambda$ vs. $A_0\eqq A$) in
  terms of throughput. Finally we are going to compare Method~1 and
  Method~2 in terms of throughput.}

\section{Solving the velocity problem} \label{Sec:Implementation}

We discuss in this section implementation and performance issues
associated with the solution of \bal{the velocity}
problem which consists of solving $A_\lambda \sfX = \sfF$ for
$\lambda \bal{\ge} 0$.

\subsection{Numerical details}

We focus our attention on approximation settings using continuous
representations of the pressure.  All the test reported in this
section and in \S\ref{Sec:solving_pressure_schur_complement},
\S\ref{Sec:Numerical_illustrations},
\S\ref{Sec:Numerical_performance_Method2} are done with continuous
Hood--Taylor Lagrange finite elements in two space dimensions on
nonuniform triangular meshes. We use Dirichlet boundary conditions
over the entire boundary of the computational domain; that is,
$\bV\eqq\bH^1_0(\Dom)$.

The parallel implementation is done with
PETSc~\citep{petsc-efficient}.  All the errors that are reported are
relative. To make sure that the results of
the paper are reproducible, we list in Table~\ref{tab:boomer_options}
all the PETSc and BoomerAMG options that are used in the paper.
In all the tests reported in the paper, the problems involving the
mass matrices $M_\bV$ and $M_Q$ are systematically solved by using
BoomerAMG with the strong threshold set to 0.1. For other matrices we
are going to set the strong threshold to either 0.1 or 0.7. The exact
value of the strong threshold will be specified in each case.  In
order to use a fixed number of V-cycles in BoomerAMG, it is required to
set the solving method to \texttt{richardson} and set the norm for
convergence testing to \texttt{none}. 

\begin{table}[h]\tt\centering
  \begin{tabular}{rl} \hline
-pc\_type & hypre\\
-pc\_hypre\_type& boomeramg \\
-pc\_hypre\_boomeramg\_strong\_threshold& 0.1 or 0.7\\
-pc\_hypre\_boomeramg\_coarsen\_type& Falgout\\
-pc\_hypre\_boomeramg\_relax\_type\_all& Chebyshev \\\hline
-ksp\_type & richardson\\
-ksp\_norm& none \\
-pc\_hypre\_boomeramg\_max\_iter& 2\\\hline
  \end{tabular}
  \caption{PETSc and BoomerAMG options used in the paper. Top: Baseline options. Bottom: additional options
 to fix the number of V-cycles. }
\label{tab:boomer_options}\vspace{-\baselineskip}
\end{table}

\bal{When solving} linear systems whose coefficients can be easily pre-computed,
like $M_Q$, $M_\bV$ or $A$, we use the iterative algebraic multigrid
solver BoomerAMG (see \cite{Henson_Yang_ANM_2002}).  Depending on the
situation, we use either a relative stopping threshold of $10^{-10}$
when accuracy is needed (for instance to invert the mass matrix $M_Q$
in \eqref{def_Schur_LA_problem} or \eqref{eq1:prop:schur_augmented_lagrange}) or we use a fixed
number of V-cycles when accuracy is not needed (typically when
inverting matrices in preconditioners).

\bal{When solving linear systems with matrices whose coefficients cannot be
easily pre-computed, we use the left preconditioned GMRES or CG methods in PETSc that only require
the action of the matrix on a vector. The relative stopping threshold for all
these situation is set to
$10^{-10}$.}

All the numerical simulations reported in the paper have been
  done in double precision on
a small cluster from the Department of Mathematics at Texas A\&M
(called ``Whistler''). This cluster in composed of 20 blades 
each composed of 2$\CROSS$Intel Xeon Gold 6130 CPUs, \@2.10GHz, with 192
GiB main memory, 10 blades composed of 2$\CROSS$Intel Xeon Gold 6226R
CPUs, \@2.10GHz, with 192 GiB main memory, and 3 blades each composed
of 2$\CROSS$Intel Xeon Gold 6226R CPUs, \@ 2.10GHz, with 768 GiB main
memory.  Each socket has 16 cores and 2 threads per core (one can run
32 MPI ranks per socket).  In principle, 66 sockets are available in
total, which make \bal{2112} MPI ranks at full load.


\subsection{Preconditioners for the velocity
  problem} \label{Sec:precond_A_plus_Lambda_D} \bal{When $\lambda>0$,}
it is very difficult to assemble the matrix coefficients of
$A_\lambda = A+ \lambda\mu B\tr M_Q^{-1} B$ \bal{because we use
  continuous finite elements to approximate the pressure} (\ie
$\polP_{k+1}/\polP_{k}$ Hood--Taylor elements, $k\ge 1$). The first
reason is that the inverse of the pressure matrix $M_Q^{-1}$ cannot be
easily computed. It can be approximated by using lumping techniques
though.  The second reason is that \bal{even if the mass is lumped,
  the sparsity pattern of the matrix is greatly increased because the
  stencil of each dof is based on two concentric layers of cells
  instead of only one}. This drastically increases the memory footprint
as well as the communications needed for parallel computing.
\bal{Therefore, we iteratively solve $A_\lambda \sfX = \sfF$ by using
  the preconditioned CG algorithm in PETSc only requiring the matrix
  action of $A_\lambda$.}

As we use Dirichlet boundary conditions, the bilinear form
$\int_\Dom\pole(\bu){:}\pole(\bv)\diff x$ is also equal to
$\int_\Dom\GRAD\bu{:}\GRAD\bv\diff x +
\int_\Dom\DIV\bu \DIV\bv\diff x$
because integrating by parts twice shows that $\int_\Dom(\GRAD\bu)\tr{:}\GRAD(\bv)\diff x =
\int_\Dom\DIV\bu\DIV\bv\diff x$ for all $\bu,\bv\in
\bH_0^1(\Dom)$. Hence defining $L_\bV$ to be the stiffness matrix
associated with the bilinear form
$\int_\Dom\GRAD\bv_h{:}\GRAD\bw_h\diff x$  and $D$ to be
the stiffness matrix associated with the grad-div bilinear form
$\int_\Dom\DIV\bv_h\DIV \bw_h\diff x$; \ie
\begin{align}
  L_{\bV,ij} &\eqq \int_\Dom \GRAD \bphi_i{:} \GRAD \bphi_j \diff x,  \qquad \forall i,j \in\calV,\\
  D_{ij} &\eqq \int_\Dom \DIV\bphi_i \DIV \bphi_j \diff x, \qquad \forall i,j \in\calV,
\end{align}
we obtain the following
equivalent representation of the matrix $A$:
\begin{equation} \label{def_of_A}
  A = \dt^{-1} M_\bV + \mu (L_\bV + D),
  \end{equation}
  where $M_\bV$ is the velocity mass matrix.  This leads use to
  consider the following three candidates to precondition the matrix
  $A_\lambda$:
\begin{align}
\tA_{\lambda,1} &\eqq \dt^{-1} M_\bV + \mu L_\bV +\mu(1+\lambda) D = A+\lambda\mu D, \label{precond_A_plus_lambda_D_1}\\
\tA_{\lambda,2} &\eqq \dt^{-1} M_\bV + \mu L_\bV +\mu D =A, \label{precond_A_plus_lambda_D_2} \\
\tA_{\lambda,3} &\eqq \dt^{-1}M_\bV + \mu L_\bV. \label{precond_A_plus_lambda_D_3} 
\end{align}
In the first preconditioner, we replace $B\tr M_Q^{-1} B$ by $D$, \bal{(see \eg \cite{Heister_Rapin_IJNMF_2013})}. In
the second preconditioner, we remove $B\tr M_Q^{-1} B$. In the third
preconditioner, we remove both $D$ and $B\tr M_Q^{-1} B$.  The stencil
for these three preconditioners is standard when using Lagrange finite
elements. Once the coefficients of the matrices
$\tA_{\lambda,1}, \tA_{\lambda,2}, \tA_{\lambda,3}$ are assembled, the
linear systems involving these matrices can be easily solved using
BoomerAMG.

\subsection{Performances of the velocity preconditioners} \label{Sec:vel_precond}

We now test the three preconditioners
\eqref{precond_A_plus_lambda_D_1}--\eqref{precond_A_plus_lambda_D_3}
with continuous Lagrange finite elements in 2D.

\subsubsection{The setting} We use the parameters
$\mu = 1 $, and $\lambda=1$.
Setting $k\eqq 16\pi$, we consider the divergence-free
vector field $\bu$ and source term $\bef$
\begin{equation}\label{divergence_free_velocity}
  \bu(\bx) =  \begin{pmatrix}\sin(k x) \sin(k y)\\\cos(k x) \cos(k y)\end{pmatrix},\qquad
  \bef(\bx) =  \frac{\bu}{\dt} -2\mu \pole(\bu).
  \end{equation}
Recalling that the velocity
shape functions are denoted $\{\bphi_i\}_{i\in\calV}$, we construct
the vector $\sfF$ with entries $\sfF_i=\int_\Dom
\bef(\bx,t)\SCAL\bphi_i(\bx)\diff x$ and solve the linear system
\begin{equation}
  A_\lambda \sfU = \sfF. \label{velocity_problem}
\end{equation}
  Setting $\bu_h\eqq \sum_{i\in\calV} \sfU_i\bphi_i$, the
relative errors in the $L^2$-norm
reported below are defined to be equal to $\|\bu-\bu_h\|_{\bL^2(\Dom)}/\|\bu\|_{\bL^2(\Dom)}$.

We perform tests with $\polP_2$ Lagrange elements
on five nonuniform meshes \bal{having respectively {237,570}, {947,714},
{3,785,730}, {15,132,674}, and {60,510,210} degrees of freedoms for the velocity}.
In order to test the weak scalability of the preconditioners, the
computations are done on those fives meshes using $2$, $8$, $32$,
$128$, and $512$ processors, respectively.

The problem \eqref{velocity_problem} is
solved with relative threshold set to $10^{-10}$ by using the CG solver in PETSc only requiring the action of $A_\lambda$.
\bal{We test two strategies to solve the linear systems
  $\tA_{\lambda,\ell} \sfX =\sfR$, \bal{$\ell\in\{1,2,3\}$}, which we
  recall is the preconditioner for \eqref{velocity_problem}.  In both
  cases we use BoomerAMG and use the entries of $\tA_{\lambda,\ell}$ to
  set the coefficients in the multigrid algorithm.  In the first
  solution strategy we let BoomerAMG iterate until the relative
  threshold $10^{-10}$ is reached. This preconditioning technique is
  denoted by $(\tA_{\lambda,\ell})\lo{th}^{-1}$. In the second strategy,}
we force BoomerAMG to use only 2 V-cycles per outer CG iteration. This
preconditioning technique is denoted by
$(\tA_{\lambda,\ell})\lo{2Vc}^{-1}$.

\subsubsection{BoomerAMG with fixed threshold}

The linear systems $\tA_{\lambda,\ell}\sfX=\sfR$, $\ell\in\{1,2,3\}$, are
solved by using the iterative algebraic multigrid solver BoomerAMG
with the options given in Table~\ref{tab:boomer_options}. We set the
strong threshold in BoomerAMG to $0.7$ for $\tA_{\lambda,1}$ and
$\tA_{\lambda,2}$
(see~\eqref{precond_A_plus_lambda_D_1}-\eqref{precond_A_plus_lambda_D_2})
and to $0.1$ for $\tA_{\lambda,3}$ (see
\eqref{precond_A_plus_lambda_D_3}). The results are shown in
Table~\ref{tab:vel_pb_exact}.  For each case, we show the number of
CG iterations that are necessary to reach the relative threshold
$10^{-10}$. We also show the $L^2$-norm of the relative error on $\bu$
and the wall-clock time in seconds for each computation.

\begin{table}[h]\scriptsize\centering
  \begin{tabular}{|c|c|c|c|c|c|c|} \hline
 \multirow{2}{*}{Precond.}& Vel. dofs & {237,570} & {947,714} & {3,785,730} & {15,132,674} &  {60,510,210} \\
                                   & Nb. Proc.     &2           & 8         & 32     & 128& 512  \\ \hline
\multirow{3}{*}{$(\tA_{\lambda,1})\lo{th}^{-1}$}& 
CG iter. & 6 & 5 & 4 & 3 & 4 \\ 
 & L2 vel. err. & 0.319E-04 & 0.217E-05 & 0.156E-06 & 0.121E-07 & 0.101E-08 \\ 
 & Times (s) & 16.5 & 25.0 & 60.0 & 83.7 & 240 \\
\hline
\multirow{3}{*}{$(\tA_{\lambda,2})\lo{th}^{-1}$}& 
CG iter. & 7 & 5 & 4 & 3 & 3 \\ 
 & L2 vel. err. & 0.319E-04 & 0.217E-05 & 0.156E-06 & 0.121E-07 & 0.101E-08 \\ 
 & Times (s) & 5.80 & 6.00 & 14.2 & 16.8 & 30.3 \\
 \hline
 \multirow{3}{*}{$(\tA_{\lambda,3})\lo{th}^{-1}$}&
CG iter. & 11 & 9 & 8 & 5 & 4 \\ 
 & L2 vel. err. & 0.319E-04 & 0.217E-05 & 0.156E-06 & 0.121E-07 & 0.101E-08 \\ 
 & Times (s) & 4.78 & 4.69 & 8.27 & 6.36 & 5.91 \\
\hline \end{tabular}
\caption{Preconditioning of
  $A_\lambda$
  (with $\lambda=1$, $\mu =1$, and $\polP_2$ Lagrange elements for the
  velocity) using $(\tA_{\lambda,1})\lo{th}^{-1}$, $(\tA_{\lambda,2})\lo{th}^{-1}$, and $(\tA_{\lambda,3})\lo{th}^{-1}$. The
  solution method for the preconditioners consists of using BoomerAMG
 and iterating on the V-cycles until the relative threshold $10^{-10}$ is reached.}%
\label{tab:vel_pb_exact}\vspace{-2\baselineskip}
\end{table}
We observe that the number of CG iterations to reach the relative
threshold $10^{-10}$ is small for the three preconditioners.  This
number converges to $1$ for the third preconditioner as the number of
grid point goes to infinity. This is because the tests are done with a
manufactured velocity field that is divergence-free (see
\S\ref{se:non_div_free}).  We also observe that the $L^2$-norm of the
relative error scales like $\calO(h^3)$ which is the expected rate of
convergence for $\polP_2$ approximation.  Notice that, although the
number of CG iterations for all the preconditioners are roughly the
same and very small, the wall-clock time for the first and second
preconditioners do not scale well.  \bal{This is due to the difficulty
  for  BoomerAMG  to construct a proper multigrid
  preconditioner for the matrix $D$ which is associated with a
  vector-valued differential operator.  Removing the contribution of
  the grad-div operator (i.e. removing the blocks above and below the
  diagonal in $\tA_\lambda$ ) to build a preconditioner
  ($(\tA_{\lambda,3})\lo{th}^{-1}$) gives a methods that scales well
  in terms of throughput.}  This test shows in passing that just
looking at the number of outer CG iterations, as sometimes done in the
literature, is not informative. The wall-clock time is also an
important factor that can help differentiate preconditioners.  In
conclusion, this series of tests clearly demonstrate that
$(\tA_{\lambda,3})\lo{th}^{-1}$ is an excellent preconditioner for the
matrix $A_\lambda$.

\subsubsection{BoomerAMG with only 2 V-cycles}

It is reported in the literature that significant CPU gains can be
obtained by degrading the solution method for the preconditioner. In
particular, some authors advocate using a fixed number of V-cycles
instead of solving $\tA_{\lambda,\ell} \sfX=\sfR$ with very high
accuracy. We now test this idea by only using two V-cycles in
BoomerAMG instead of letting the algorithm run until reaching the
assigned tolerance.  The results obtained with this preconditioning
strategy are shown in Table~\ref{tab:vel_pb_2fvc}.

\begin{table}[h]\scriptsize\centering
  \begin{tabular}{|c|c|c|c|c|c|c|} \hline
 \multirow{2}{*}{Precond.}& Vel. dofs & {237,570} & {947,714} & {3,785,730} & {15,132,674} &  {60,510,210} \\ 
                         & Nb. Proc.     &2           & 8         & 32     & 128& 512  \\ \hline
\multirow{3}{*}{$(\tA_{\lambda,1})\lo{2Vc}^{-1}$}& 
CG iter. & 29 & 65 & 738 & 1030 & 2825 \\ 
 & L2 vel. err. & 0.319E-04 & 0.217E-05 & 0.805E+05 & 0.111E+06 & 0.105E+06 \\ 
 & Times (s) & 5.15 & 12.7 & 243 & 356 & 1049 \\
 \hline
\multirow{3}{*}{$(\tA_{\lambda,2})\lo{2Vc}^{-1}$}& 
CG iter. & 41 & 39 & 114 & 148 & 285 \\ 
 & L2 vel. err. & 0.319E-04 & 0.217E-05 & 0.156E-06 & 0.121E-07 & 0.101E-08 \\ 
 & Times (s) & 6.16 & 7.23 & 34.4 & 48.5 & 96.9 \\
 \hline
 \multirow{3}{*}{$(\tA_{\lambda,3})\lo{2Vc}^{-1}$}&
CG iter. &18 & 20 & 21 & 23 & 24 \\ 
 & L2 vel. err. & 0.319E-04 & 0.217E-05 & 0.156E-06 & 0.121E-07 & 0.101E-08 \\ 
 & Times (s) & 2.36 & 2.91 & 5.47 & 6.76 & 7.61  \\
\hline \end{tabular}
\caption{Preconditioning of $A_\lambda$ (with $\lambda=1$, $\mu =1$,
  and $\polP_2$ Lagrange elements for the velocity) using
  $(\tA_{\lambda,1})\lo{2Vc}^{-1}$, $(\tA_{\lambda,2})\lo{2Vc}^{-1}$, and $(\tA_{\lambda,3})\lo{2Vc}^{-1}$.  The
  solution method for the preconditioners consists of using BoomerAMG
  with 2 V-cycles.}%
\label{tab:vel_pb_2fvc}\vspace{-2\baselineskip}
\end{table}

We observe that using 2 V-cycles to accelerate the performance of the
preconditioners significantly increases the number of outer CG
iterations.  For the first two preconditioners,
$(\tA_{\lambda,1})\lo{2Vc}^{-1}$, $(\tA_{\lambda,2})\lo{2Vc}^{-1}$, the
number of CG iterations even diverges when refining the mesh. For
reasons we do not fully understand, this behavior seems to be related
to our using a divergence-free manufactured solution (see
\S\ref{se:non_div_free} for other details in this directions).  We
observe however that $(\tA_{\lambda,3})\lo{2Vc}^{-1}$ behaves properly,
and even though the number of CG iterations is larger than when
using $(\tA_{\lambda,3})\lo{th}^{-1}$, the actual CPU times is similar to
what is reported in Table~\ref{tab:vel_pb_exact}.

This series of tests confirms again that $(\tA_{\lambda,3})\lo{2Vc}^{-1}$
is an excellent preconditioner and suggests that using only 2 V-cycles
in BoomerAMG is a viable option for this preconditioner, whereas it
does not seem to be the case of
$(\tA_{\lambda,1})\lo{2Vc}^{-1}$ and
$(\tA_{\lambda,2})\lo{2Vc}^{-1}$.  We will see though in
\S\ref{Sec:Numerical_performance_Method2} that
$(\tA_{\lambda,2})\lo{2Vc}^{-1}$ is a good approximate preconditioner of
$A_\lambda$ when used in the context of the solution Method~2 (described
\S\ref{Sec:Method_2}).

\subsection{Tests with non divergence-free solution} \label{se:non_div_free}

We now confirm the conjecture made in \S\ref{Sec:vel_precond} that the
number of CG iterations to solve the problem
\eqref{velocity_problem} converges to 1 as the number of grid point
goes to infinity because the solution defined
in~\eqref{divergence_free_velocity} is divergence-free.

In this section we test the third preconditioner, $\tA_{\lambda,3}$,
with a right-hand side $\sfF$ constructed with a velocity field that
is not divergence-free. More precisely, we take
\begin{equation}\label{non_divergence_free_velocity}
  \bu(\bx) =  \begin{pmatrix}\sin(2k x) \sin(k y)\\\cos(k x) \cos(k y)\end{pmatrix},\qquad
  \bef(\bx) =  \frac{\bu}{\dt} -2\mu \pole(\bu) -\lambda\mu\GRAD\DIV\bu,
\end{equation}
with $k=16\pi$, $\lambda=1$, $\mu=1$, and $\dt=N^{-\frac12}$.  Here
again we solve $A_\lambda \sfU = \sfF$ and the entries of the vector
$\sfF$ are given by
$\sfF_i=\int_\Dom \bef(\bx,t)\SCAL\bphi_i(\bx)\diff x$. The tests are
done with $\polP_2$ and $\polP_3$ finite elements. For the $\polP_3$
elements, the tests are done with meshes composed of
\bal{{133,874}, {533,570}, {2,130,434}, {8,514,050}, {34,040,834} dofs} 
respectively. Recall that the number of degrees of freedom is twice
the number of grid points as we are solving a vector-valued problem in
dimension two.

\begin{table}[h]\scriptsize\centering
   \begin{tabular}{|c|c|c|c|c|c|c|c|} \hline
 F.E. & Precond. &  Nb. Proc.     &2           & 8         & 32     & 128& 512  \\ \hline
 \multirow{7}{*}{ $\polP_2$ } &  & Vel. dofs& {237,570} & {947,714} & {3,785,730} & {15,132,674} &  {60,510,210} \\
  \cline{2-8}
 & \multirow{4}{*}{$(\tA_{\lambda,3})\lo{th}^{-1}$}
 & CG iter. & 34 & 28 & 25 & 20 & 18 \\ 
 && L2 vel. err. & 0.279E-03 & 0.185E-04 & 0.127E-05 & 0.923E-07 & 0.712E-08 \\ 
 && Times (s) & 13.2 & 12.8 & 23.2 & 19.7 & 18.9  \\
 \cline{2-8}
 &\multirow{3}{*}{$(\tA_{\lambda,3})\lo{2Vc}^{-1}$}& 
CG iter. & 40 & 32 & 27 & 20 & 18 \\ 
 && L2 vel. err. & 0.279E-03 & 0.185E-04 & 0.127E-05 & 0.923E-07 & 0.712E-08 \\ 
 && Times (s) & 5.10 & 4.52 & 6.81 & 5.76 & 5.63  \\ \hline
 \multirow{7}{*}{ $\polP_3$ } & & Vel. dofs & \phantom{0}{133,874}&{533,570}&{2,130,434}&{8,514,050}&{34,040,834} \\
  \cline{2-8}
  & \multirow{3}{*}{$(\tA_{\lambda,3})\lo{th}^{-1}$}&
CG iter. & 35 & 31 & 26 & 23 & 18 \\ 
 && L2 vel. err. & 0.512E-03 & 0.405E-04 & 0.272E-05 & 0.174E-06 & 0.109E-07 \\ 
 && Times (s) & 12.4 & 12.8 & 19.8 & 18.9 & 16.3   \\
 \cline{2-8}
  &\multirow{3}{*}{$(\tA_{\lambda,3})\lo{2Vc}^{-1}$}&
CG iter. & 68 & 49 & 36 & 24 & 20 \\ 
 && L2 vel. err. & 0.512E-03 & 0.405E-04 & 0.272E-05 & 0.174E-06 & 0.109E-07 \\ 
 && Times (s) & 7.56 & 6.09 & 7.46 & 5.84 & 5.47 \\
\hline \end{tabular}
\caption{Non divergence-free solution. Preconditioning of $A_\lambda$
  with $(\tA_{\lambda,3})\lo{th}^{-1}$ or $(\tA_{\lambda,3})\lo{2Vc}^{-1}$, (with $\dt=N^{-\frac12}$, $\lambda=1$,
  $\mu=1$) using BoomerAMG with strong threshold $0.1$. Top: $\polP_2$ 
  Lagrange elements. Bottom: $\polP_3$ Lagrange
  elements.}
  \label{tab:exact_vs_2fvc_ndv_sol}%
  \vspace{-\baselineskip}
\end{table}

The results are shown in Table~\ref{tab:exact_vs_2fvc_ndv_sol}. 
\bal{Contrary to} what is observed in
Tables~\ref{tab:vel_pb_exact}, the number of
CG iterations for $(\tA_{\lambda,3})\lo{th}^{-1}$ \bal{now seems to converge to
a number larger than 1.}  The $L^2$-norm of the
relative error scales like $\calO(h^3)$ for the $\polP_2$
approximation and like $\calO(h^4)$ for the $\polP_3$ approximation
until the $10^{-10}$ threshold is reached, as expected.

This time, however, fixing the number of V-cycles in BoomerAMG does
not significantly increases the number of CG iterations that are
necessary to reach the assigned threshold on fine meshes. Looking at
lines lines 3\&6 and lines 10\&13 in
Table~\ref{tab:exact_vs_2fvc_ndv_sol}, we observe that the number of
CG iterations only increases by one or two units (contrary to what is
observed when comparing line 9 in Table~\ref{tab:vel_pb_exact} and
line 9 in Table~\ref{tab:vel_pb_2fvc}).  As a result, fixing the
number of V-cycles in BoomerAMG speeds-up the computation by a factor
3 (compare lines 5\&8 and lines 12\&15 in
Table~\ref{tab:exact_vs_2fvc_ndv_sol}).  Notice that the performance of the 
\bal{preconditioner $(\tA_{\lambda,3})^{-1}\lo{2VC}$}
on  the
divergence-free solution and the non divergence-free solution are
similar (compare line 8 in Table~\ref{tab:exact_vs_2fvc_ndv_sol}
with line 11 in Table~\ref{tab:vel_pb_2fvc}).

Finally, we observe that the weak scaling is satisfactory for both the
$\polP_2$ and $\polP_3$ approximations.


\subsection{Conclusions for \S\ref{Sec:Implementation}} We draw two important conclusions from the
series of tests done in section \S\ref{Sec:Implementation}. The first
one is that although the grad-div operator has very interesting
properties emphasized in the literature
(\cite{Benzi_Olshanksii_SISC_2006,Olshanskii_Lube_Heister_Lowe_CMAME_2009},
\cite{Guermond_Minev_SISC_2015,deFrutos_Garcia_John_Novo_JSC_2016}),
 efficiently solving problems requiring the solution of linear systems
involving the matrix $D$ or variations thereof is difficult (see \eg
\cite{Heister_Rapin_IJNMF_2013,Jenkins_John_Linke_Rebholz_ACM_2014,Fiordilino_Layton_Rong_CMAME_2018}).
The best preconditioner for this problem is $\tA_{\lambda,3}$, which
consists of dropping the coupling block $D$.
The second conclusion is
that a speed-up factor between 2 and 3 can be gained by requesting that BoomerAMG
only performs 2 V-cycles at each iteration of the velocity problem, as done in $(\tA_{\lambda,3})\lo{2Vc}^{-1}$,
instead of reaching a fixed relative threshold.

\section{Solving the pressure Schur complement} \label{Sec:solving_pressure_schur_complement}

In this section we consider various preconditioners of the pressure
Schur complement $S_\lambda$ defined in \eqref{def_Schur_LA_problem}.
These preconditioners are tested in
\S\ref{Sec:Numerical_illustrations} and
\S\ref{Sec:Numerical_performance_Method2}.

\subsection{Preconditioning $S_\lambda$}

Ideas to precondition $S_{\lambda}$ can be found by investigating how
$S_{\lambda}$ degenerates in the limit $\mu\dt/h^2\to \infty$
and in the limit $\mu\dt/h^2\to 0$; see \eg \cite{Cahouet_Chabard_1988} and
\cite{Bramble_Pasciak_1997}.

In the limit $\mu\dt/h^2\to 0$, we have  $S_{\lambda}\to  \dt B (M_\bV)^{-1}
B\tr$. The inverse of $S_{\lambda}$ in this limit is then $\dt^{-1} (B (M_\bV)^{-1}
B\tr)^{-1}$.

In the limit $\mu\dt/h^2 \to \infty$, we have
$S_{\lambda}\to \mu^{-1} B(E_\bV + \lambda B\tr M_Q^{-1} B)^{-1}
B\tr$, where we recall that the matrix $E_\bV$ is defined in
\eqref{def_of_EV}. As $E_\bV$ and $B\tr M_Q^{-1} B$ have similar
spectra, it seems natural to approximate $S_{\lambda}$ by
$\mu^{-1} (1+\lambda)^{-1} B(B\tr M_Q^{-1} B)^{-1} B\tr$ in the limit
$\mu\dt/h^2 \to \infty$. Hence we
posit that the inverse of $S_{\lambda}$ approximately behaves like
$\mu (1+\lambda)M_Q^{-1}$ when $\mu\dt/h^2 \to \infty$.

In conclusion, the first preconditioner for $S_{\lambda}$ that we consider in the paper
is
\begin{equation}
  S_{\lambda}^{-1} \sim  C_{\lambda}\eqq \dt^{-1} (B (M_\bV)^{-1} B\tr)^{-1} + \mu (1+\lambda)M_Q^{-1}.
\label{AL_preconditioner}
\end{equation}
In light of the identity \eqref{eq1:prop:schur_augmented_lagrange},
the above expression is consistent with the approximation of
$S_0^{-1}$ by $\dt^{-1} (B (M_\bV)^{-1} B\tr)^{-1} + \mu M_Q^{-1}$,
which is the preconditioner of the Schur complement $S_0$ considered
in \cite[Eq.~(44)]{Cahouet_Chabard_1988}.

\bal{As the matrix $B M_\bV^{-1} B\tr$}
is only used as part of the preconditioner for $S_{\lambda}$,
replacing the mass matrix $M_\bV$ by a lumped one could save \bal{some
computation} time. To explore this idea, we  propose to replace $M_\bV$ by a diagonal
matrix $\Lambda_\bV$ with entries defined by%
\begin{equation}%
\Lambda_{\bV,ij} \eqq \delta_{ij}  \int_\Dom|\varphi_i(\bx)|\diff x,\qquad \forall i,j\in\calV.
  \end{equation}
Taking the absolute value of the shape function is important here as
$\polP_2$ and $\polP_3$ velocity shape functions are not uniformly positive.

In \citep{Cahouet_Chabard_1988}, the authors also propose using the
discrete pressure Laplacian instead of $B(M_\bV)^{-1} B\tr$, \bal{\ie they propose the preconditioner
$  S_{\lambda}^{-1} \sim \dt^{-1} (L_Q)^{-1} + \mu (1+\lambda)M_Q^{-1}$ where
$L_Q$ is the matrix of the pressure Laplacian with entries}
\begin{equation}
L_{Q,kl}\eqq  \int_\Dom \GRAD \psi_k\SCAL\GRAD \psi_l \diff x\qquad \forall k,l,\in\calQ. \label{pressure_Laplacian}
\end{equation}

In conclusion, in the rest of this paper we are going to compare the
following three preconditioners for $S_{\lambda}$:
\begin{subequations} \label{Schur_complement_preconditioners}
  \begin{align}
    C_\lambda&\eqq \dt^{-1} (B\ M_\bV^{-1} B\tr)^{-1} + \mu (1+\lambda) M_Q^{-1},\\
    C_\lambda^\Lambda&\eqq\dt^{-1} (B\ \Lambda_\bV^{-1} B\tr)^{-1} + \mu (1+\lambda) M_Q^{-1}, \\ 
    C_\lambda^\Delta&\eqq\dt^{-1} (L_Q)^{-1} + \mu (1+\lambda) M_Q^{-1}.\label{Slambda_Lap_preconditioner}
  \end{align}
\end{subequations}

In the next two sections we discuss practical details regarding the
implementation of the above preconditioners using PETSc and BoomerAMG.

\subsection{Solving $B M_V^{-1}B\tr \sfX = \sfY$} \label{Sec:BMvmBt}


The coefficients of the matrix $B M_\bV^{-1} B\tr$ cannot be easily
computed when using continuous elements for the pressure. Even if the
velocity mass matrix is lumped, the stencil of the resulting
simplified matrix \bal{is based on two concentric layers of cells
  around each dof.} Hence, the only practical way to solve the linear
system $B M_\bV^{-1} B\tr \sfX = \sfY$ is to use a preconditioned
iterative method.  We use the PETSc version of GMRES only requiring
the matrix action.

We use BoomerAMG with strong threshold equal to 0.1 to estimate
$M_\bV^{-1} \sfZ$ (\ie to solve $M_\bV \sfW = \sfZ$). Henceforth, when
we write $(M_{\bV})\lo{th}^{-1} \sfZ$, we mean that the linear system
$M_\bV \sfW = \sfZ$ is solved with BoomerAMG with relative threshold
$10^{-10}$. As $B M_\bV^{-1} B\tr$ is just a preconditioner, we can
also replace the solution of the linear system $M_\bV \sfW = \sfZ$ by
just using 2 V-cycles in BoomerAMG. We are going to use the notation
$B (M_{\bV})\lo{2Vc}^{-1} B\tr$ to refer to this other method.

A natural preconditioner for the matrices $B(M_{\bV})\lo{th}^{-1} B\tr$ and
$B (M_{\bV})\lo{2Vc}^{-1} B\tr$ is the pressure Laplacian matrix
$L_Q$ defined in~\eqref{pressure_Laplacian}.  However, this matrix is
singular when Dirichlet boundary conditions are enforced on the
velocity over the entire boundary of the domain $\Dom$. As we only
want to construct a preconditioner, one can consider instead the
matrix with entries $\epsilon M_Q+L_Q$ with $\epsilon>0$. We henceforth set
$\epsilon=1$.

In the rest of this section we
evaluate the performance of the preconditioner $\epsilon M_Q+L_Q$ using BoomerAMG (with
strong threshold equal to 0.1).
We do two series of tests. In the first series, we
precondition the matrix $B (M_\bV)\lo{th}^{-1} B\tr$
with $\epsilon M_Q+L_Q$, and
we solve the linear system $(\epsilon M_Q+L_Q)\sfW =\sfZ$ with
BoomerAMG by iterating until the relative threshold $10^{-10}$ is reached. We use the
notation$(\epsilon M_Q+L_Q)\lo{th}^{-1}$ for this preconditioner.
In the second series of tests, we precondition the matrix
$B (M_{\bV})\lo{2Vc}^{-1} B\tr$ with 2 V-cycles of BoomerAMG using
the coefficients of the matrix $\epsilon M_Q+L_Q$. We use the notation
$(\epsilon M_Q+L_Q)\lo{2Vc}^{-1}$ for this preconditioner.

The tests are done with continuous Lagrange finite elements in two
space dimensions.  We consider mixed $\polP_2/\polP_1$ and
$\polP_3/\polP_2$ elements.
For each mesh we build the
right-hand $\sfY\eqq B (M_\bV)\lo{th}^{-1} B\tr \sfX$ (or
$\sfY\eqq B (M_{\bV})\lo{2Vc}^{-1} B\tr \sfX$) by applying the
matrix $B (M_\bV)\lo{th}^{-1} B\tr$ (or
$\sfY\eqq B (M_{\bV})\lo{2Vc}^{-1} B\tr$) to the vector $\sfX$
whose entries are the values of the following pressure field at the
Lagrange nodes of the pressure mesh:
\begin{equation}
  p(x,y) \eqq \sin(16\pi(x-y)). \label{exact_pressure}
\end{equation}

\begin{table}[h]\scriptsize\centering\addtolength{\tabcolsep}{-2.2pt}
   \begin{tabular}{|c|c|c|c|c|c|c|c|} \hline
 \multirow{2}{*}{F.E.} &  \multirow{2}{*}{Method.} & Press. dofs & 29,857 & 11,8785 & 47,3857 & 189,2865 &  7,566,337 \\ 
                                    & & Nb. Proc.     &2           & 8         & 32     & 128& 512  \\ \hline
 \multirow{6}{*}{ $\polP_1$ } & \multirow{3}{*}{$(\epsilon M_Q+L_Q)\lo{th}^{-1}$}&
GMRES iter. & 9 & 9 & 8 & 8 & 8 \\
&&L1 press. err. & 0.478E-10 & 0.239E-10 & 0.470E-10 & 0.643E-10 & 0.398E-10 \\
&&Time (s) & 2.532 & 2.911 & 5.057 & 6.360 & 6.439  \\
 \cline{2-8}
 &\multirow{3}{*}{$(\epsilon M_Q+L_Q)\lo{2Vc}^{-1}$}&
GMRES iter. & 9 & 9 & 8 & 8 & 8 \\
&&L1 press. err. & 0.475E-10 & 0.235E-10 & 0.470E-10 & 0.641E-10 & 0.399E-10  \\
&&Time (s) & 0.843 & 0.983 & 1.556 & 2.356 & 2.435    \\
 \hline
 \multirow{6}{*}{ $\polP_2$ } & \multirow{3}{*}{$(\epsilon M_Q+L_Q)\lo{th}^{-1}$}&
GMRES iter. &  11 & 10 & 10 & 9 & 8 \\
&&L1 press. err. & 0.633E-10 & 0.323E-10 & 0.472E-11 & 0.519E-11 & 0.111E-10  \\
&&Time (s) & 2.835 & 3.001 & 5.469 & 6.712 & 6.545  \\
 \cline{2-8}
 &\multirow{3}{*}{$(\epsilon M_Q+L_Q)\lo{2Vc}^{-1}$}&
GMRES iter. &  12 & 11 & 11 & 10 & 10   \\
&&L1 press. err. & 0.279E-10 & 0.673E-10 & 0.189E-10 & 0.434E-10 & 0.327E-10 \\
&&Time (s) & 0.782 & 0.823 & 1.281 & 2.432 & 2.222 \\
\hline \end{tabular}%
\caption{Tests for the preconditioning of $B (M_\bV)\lo{th}^{-1} B\tr$ by
  $(\epsilon M_Q+L_Q)\lo{th}^{-1}$ and the preconditioning of $B (M_\bV)\lo{2Vc}^{-1} B\tr$ by
  $(\epsilon M_Q+L_Q)\lo{2Vc}^{-1}$ with $\polP_1$ and
  $\polP_2$ Lagrange elements, using BoomerAMG, strong threshold
  0.1.}%
\label{tab:BmBt_exact_2fvc}%
\end{table}%

We test five (mixed) meshes.  The number of pressure grid points on
each of these meshes are {29,857}, {118,785}, {473,857}, {1,892,865},
and {7,566,337}.  The results of these tests are shown in
Table~\ref{tab:BmBt_exact_2fvc}. We show the number of GMRES
iterations, the $L^1$-norm of the relative error on the pressure when
the GMRES threshold is reached, and the wall-clock time. As the GMRES
threshold is $10^{-10}$, \bal{we observe that the errors are of order
$10^{-10}$} (these are not approximation errors since there is nothing
to approximate).

This test shows that $(\epsilon M_Q + L_Q)\lo{th}^{-1}$ and
$(\epsilon M_Q + L_Q)\lo{2Vc}^{-1}$ are excellent preconditioners of
$B (M_\bV)\lo{th}^{-1} B\tr$ and $B (M_\bV)\lo{2Vc}^{-1} B\tr$,
respectively, as the number of GMRES iterations is small and decreases
as the mesh is refined.

\subsection{Solving $B \Lambda_V^{-1}B\tr \sfX = \sfY$}  \label{Sec:B_LambdaV_Bt}
We now investigate solution methods to solve
$B \Lambda_\bV^{-1} B\tr\sfP=\sfF$. Note that this time the inverse of
$\Lambda_\bV$ can be calculated exactly. Here again we can
precondition the system with either
$(\epsilon M_Q + L_Q)\lo{th}^{-1}$ or
$(\epsilon M_Q + L_Q)\lo{2Vc}^{-1}$.

\begin{table}[h]\scriptsize\centering\addtolength{\tabcolsep}{-2.2pt}
   \begin{tabular}{|c|c|c|c|c|c|c|c|} \hline
 \multirow{2}{*}{F.E.} &  \multirow{2}{*}{Method.} & Press. dofs & 29857 & 118785 & 473857 & 1892865 &  7566337 \\ 
                                    & & Nb. Proc.     &2           & 8         & 32     & 128& 512  \\ \hline
 \multirow{6}{*}{ $\polP_1$ } & \multirow{3}{*}{$(\epsilon M_Q+L_Q)\lo{th}^{-1}$,}&
GMRES iter. & 15 & 15 & 14 & 14 & 14  \\
&&L1 press. err. & 0.504E-10 & 0.128E-09 & 0.534E-10 & 0.663E-09 & 0.125E-08  \\
&&Time (s) & 0.852 & 1.01 & 1.48 & 2.61 & 2.44  \\
 \cline{2-8}
 &\multirow{3}{*}{$(\epsilon M_Q+L_Q)\lo{2Vc}^{-1}$}&
GMRES iter. & 15 & 15 & 14 & 15 & 15 \\
&&L1 press. err. & 0.531E-10 & 0.325E-10 & 0.582E-10 & 0.960E-10 & 0.233E-09\\
&&Time (s) & 0.648 & 0.733 & 0.981 & 1.98 & 1.97  \\
 \hline
 \multirow{6}{*}{ $\polP_2$ } & \multirow{3}{*}{$(\epsilon M_Q+L_Q)\lo{th}^{-1}$,}&
GMRES iter. & 18 & 18 & 17 & 17 & 16  \\
&&L1 press. err. & 0.710E-10 & 0.807E-10 & 0.948E-10 & 0.137E-09 & 0.545E-09 \\
&&Time (s) & 1.15 & 1.39 & 2.11 & 3.16 & 3.99 \\
 \cline{2-8}
 &\multirow{3}{*}{$(\epsilon M_Q+L_Q)\lo{2Vc}^{-1}$}&
GMRES iter. & 18 & 18 & 17 & 17 & 17   \\ 
&&L1 press. err. & 0.757E-10 & 0.603E-10 & 0.711E-10 & 0.132E-09 & 0.127E-09 \\
&&Time (s) & 0.830 & 0.928 & 1.23 & 2.30 & 2.61 \\
\hline \end{tabular}
\caption{Tests for the preconditioning of $B \Lambda_\bV^{-1} B\tr$ with $(\epsilon M_Q+L_Q)\lo{th}^{-1}$
  and $(\epsilon M_Q+L_Q)\lo{2Vc}^{-1}$
with  $\polP_1$ and $\polP_2$ Lagrange elements, using BoomerAMG, strong threshold 0.1.}
  \label{tab:BmBt_exact_2fvc_lumped}%
  \vspace{-\baselineskip}
\end{table}

We perform the same tests as in \S\ref{Sec:BMvmBt}. We solve the
problem $B\Lambda_\bV^{-1} B\tr \sfX =\sfY$ where the right-hand side
is computed by setting the entries of $\sfX$ to be the values of the
scalar field given in \eqref{exact_pressure} at the Lagrange nodes of
the pressure mesh.  The results of this series of tests are shown in
Table~\ref{tab:BmBt_exact_2fvc_lumped}. We observe that even though
using $\Lambda_\bV$ instead of $M_\bV$ takes more GMRES iterations to
achieve convergence, the wall-clock time is shorter by 30\%. Once
again the tests are purely algebraic; as a result, the errors are
proportional to the GMRES threshold.
\section{Numerical illustration of the performance of Method~1}
\label{Sec:Numerical_illustrations}

In this section we focus on the solution Method 1 defined
in~\eqref{Schur_LA_problem} and compare the performance of the
preconditioners $C_\lambda$, $C_\lambda^\Lambda$, and
$C_\lambda^\Delta$ defined \eqref{Schur_complement_preconditioners}.

\subsection{Notation} \label{Sec:notation} To be precise and make sure that the tests
reported in the paper are reproducible, we precisely define the
preconditioners we are going to use for the pressure Schur
complement. Using the symbols ``$\mathsf{th}$'' and ``$\mathsf{2Vc}$''
introduced in the previous sections, for all
$\sfa \in \{\mathsf{th},\mathsf{2Vc}\}$, we set
\begin{subequations} \label{Schur_complement_preconditioners_num}
  \begin{align}
   C_{\lambda,\sfa}&\eqq \mu (1+\lambda) (M_Q)_{\sfa}^{-1}+ \dt^{-1} (B (M_\bV)_{\sfa}^{-1} B\tr)_{\sfa}^{-1},\\
    C_{\lambda,\sfa}^\Lambda&\eqq \mu (1+\lambda) (M_Q)_{\sfa}^{-1} + \dt^{-1} (B \Lambda_\bV^{-1} B\tr)_{\sfa}^{-1}, \\
    C_{\lambda,\sfa}^\Delta&\eqq \mu (1+\lambda) (M_Q)_{\sfa}^{-1} + \dt^{-1} (L_Q)_{\sfa}^{-1}.
                      \label{Slambda_Lap_preconditioner_num}
  \end{align}
\end{subequations}

We recall that the notation $\sfW=(M_Q)_{\sfa}^{-1}\sfZ$ means that we
solve the problem $M_Q \sfW = \sfZ$ with BoomerAMG by iterating until
the relative residual reaches the threshold $10^{-10}$ if
$\sfa=\mathsf{th}$ or by using 2 V-cycles if $\sfa=\mathsf{2Vc}$.  The
same convention holds for $\sfW=(L_Q)_{\sfa}^{-1}\sfZ$ and
$\sfW=(M_\sfV)_{\sfa}^{-1}\sfZ$.

\bal{The notation $\sfW=(B (M_\bV)_{\sfa}^{-1} B\tr)_{\sfa}^{-1}\sfZ$
  means that the problem
  $B (M_\bV)_{\sfa}^{-1} B\tr \sfW= \sfZ$ is solved with GMRES using the matrix
  action and preconditioned by $(\epsilon M_Q+L_Q)\lo{\sfa}^{-1}$. We use the same definition for
 $\sfW=(B \Lambda_\bV^{-1} B\tr)_{\sfa}^{-1}\sfZ$.
We refer the reader to \S\ref{Sec:BMvmBt} and \S\ref{Sec:B_LambdaV_Bt} where we
explain how the problems $(B (M_\bV)_{\sfa}^{-1} B\tr)W=\sfZ$ and
$(B \Lambda_\bV^{-1} B\tr)_{\sfa}W=\sfZ$ are solved. }

\subsection{Numerical details} \label{Sec:Numerical_details}

We recall that to be representative of situations corresponding to the
approximation of the time-dependent Navier-Stokes equations, where the
time step decreases like the mesh size due to the nonlinearities being
made explicit in time, we set $\tau=N^{-\frac12}$, where $N$ is the
total number of grid points for the velocity.

Method 1 consists of iterating on the pressure Schur complement,
$S_\lambda$, using GMRES until the relative threshold on the residual
is less than or equal to $10^{-10}$. The
following operations are done at each GMRES iterations:
\begin{itemize}
\item One matrix-vector multiplication with the Schur complement matrix
  $S_{\lambda}$. The only nontrivial step in this operation consists of solving one
  velocity problem $A_\lambda \sfX = \sfY$ where $A_\lambda$ is defined in
  \eqref{def_of_Alambda}. This is done with CG only using the matrix action. The
        relative threshold is $10^{-10}$. The preconditioner is
        $(\tA_{\lambda,3})\lo{2Vc}^{-1}$. This  requires the following
          non-trivial operations:
        \begin{itemize}
        \item If $\lambda>0$, the matrix-vector multiplications by $A_\lambda$ require
          solving a pressure mass problem. This is done using $(M_Q)\lo{th}^{-1}$.
        \item The application of the preconditioner
          $(\tA_{\lambda,3})\lo{2Vc}^{-1}$.
        \end{itemize}
   
      \item One application of the preconditioner, \ie
        $C_{\lambda,\sfa,\sfb}$, $C^\Lambda_{\lambda,\sfa}$, or
        $C^ \Delta_{\lambda,\sfa,\sfb}$. This requires:
        \begin{itemize}
        \item Solving one pressure mass problem $(M_Q)\lo{a}^{-1}\sfZ$.
        \item Solving one problem involving
          \bal{$(B(M_{\bV})\lo{a}^{-1}B^T)^{-1}\lo{a}$ or
          $(B\Lambda_V^{-1}B^T)^{-1}\lo{a}$.} This is done using GMRES with
          rel. thr. $10^{-10}$ and preconditioning the system with
          $(\epsilon M_Q + L_V)\lo{a}^{-1}$. This requires the
          following non-trivial operations:
            \begin{itemize}
            \item Solving a velocity mass problem $(L_\bV)\lo{a}^{-1}\sfZ$ if one uses \bal{$B(M_{\bV})\lo{a}^{-1}B^T$.}
           \item The application of the preconditioner
             $(\epsilon M_Q + L_V)_{\textup{th}}^{-1}$.
            \end{itemize}
        \end{itemize}
\end{itemize}

All the simulations are done with mixed
$\polP_2/\polP_1$ continuous finite elements.
We test five meshes with
the following velocity/pressure grid point counts: Mesh~1
(118,785/29,857); Mesh~2 (473,857/118,785); Mesh~3 (1,892,865/473,857);
Mesh~4 (7,566,337/1,892,865); Mesh~5 (30,255,105/7,566,337).
As the velocity field is two dimensional,
the total number of degrees of freedom is for each mesh: Mesh~1
(267,427); Mesh~2 (1,066,499); Mesh~3 (4,259,587); Mesh~4
(17,025,539); Mesh~5 (68,076,547).
Computation done on $\polP_3/\polP_2$  finite elements gave similar results as the ones shown here
and are therefore not reported for brevity.

In order to properly compare the various methods,
we estimate the \bal{throughput per second that is achieved with the relative threshold $10^{-10}$.
This is the ratio consisting of dividing the total number of degrees of freedom by
the wall-clock time (in second) multiplied by the number of processors:}
\begin{equation}
  \text{TPS} \eqq \frac{
  \text{nb. vel.+press. degrees of freedom}}{\text{wall-clock time(s)}\CROSS \text{nb. processor}}.
\end{equation}
The inverse of $\text{TPS}$ is the computational time (s) spent per
degrees of freedom to achieve the desired relative threshold on the
residual (\underline{$10^{-10}$ in our case}).

To avoid making a compilation of dozens of tables, we have chosen to
only focus our attention on the throughput per second.  We thus do not
discuss the convergence of the errors on the velocity and pressure but
we have verified that they scales properly with the mesh size.

\subsection{Tests without augmented Lagrangian}

We start by illustrating the method proposed in
\cite[Eq.~(44)]{Cahouet_Chabard_1988}.  This method does not use the
augmented Lagrangian; that is, the method they propose consists of
setting $\lambda=0$ in \eqref{Schur_LA_problem} and
\eqref{AL_preconditioner}. The authors consider the following
preconditioner for $S_{0}$:
\begin{subequations}
  \begin{align}
    C_0& \eqq   \mu  M_Q^{-1}+ \dt^{-1} (B\ M_\bV^{-1} B\tr)^{-1}.
  \end{align}
But, based on the results shown in \S\ref{Sec:B_LambdaV_Bt}, we will also test the preconditioner,
\begin{equation}
  C_0^\Lambda \eqq  \mu  M_Q^{-1} + \dt^{-1} (B\ \Lambda_\bV^{-1} B\tr)^{-1} .
  \end{equation}
\end{subequations}
We are going to test two versions for each of these preconditioners:
\begin{subequations}\label{Cahouet_Chabard_precondiotioners}
  \begin{align}
(C_0)\lo{th} &\eqq \mu  (M_Q)\lo{th}^{-1}
    + \dt^{-1} (B\ (M_\bV)\lo{th}^{-1} B\tr)\lo{th}^{-1} ,\\
(C_0)\lo{2Vc} &\eqq \mu  (M_Q)\lo{2Vc}^{-1}
    + \dt^{-1} (B\ (M_\bV)\lo{2Vc}^{-1} B\tr)\lo{2Vc}^{-1} , \\
  (C_0^\Lambda)\lo{th} &\eqq \mu  (M_Q)\lo{th}^{-1}
    + \dt^{-1} (B\ \Lambda_\bV^{-1} B\tr)\lo{th}^{-1} ,\\
(C_0^\Lambda)\lo{2Vc} &\eqq \mu  (M_Q)\lo{2Vc}^{-1}
    + \dt^{-1} (B\ \Lambda_\bV^{-1} B\tr)\lo{2Vc}^{-1}.
 \end{align}
\end{subequations}

\begin{figure}[h]
  \centering
  \includegraphics[height=0.24\textwidth,trim={0 10 44 0},clip]{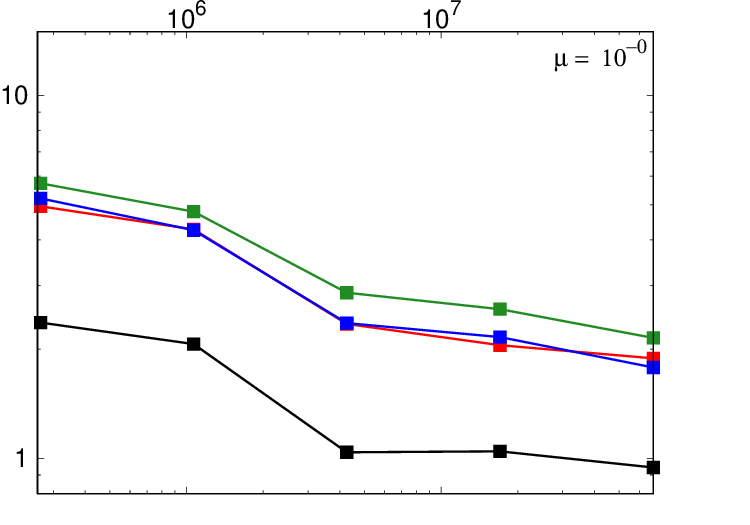}
  \includegraphics[height=0.24\textwidth,trim={16 10 44 0},clip]{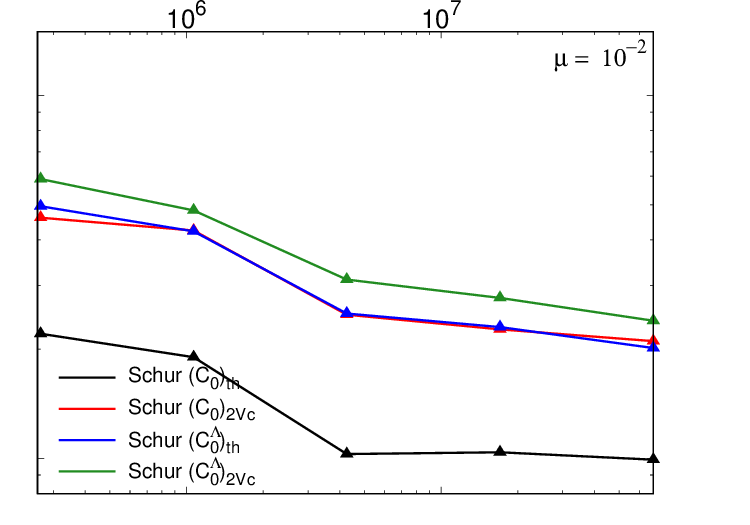}
  \includegraphics[height=0.24\textwidth,trim={16 10 0 0},clip]{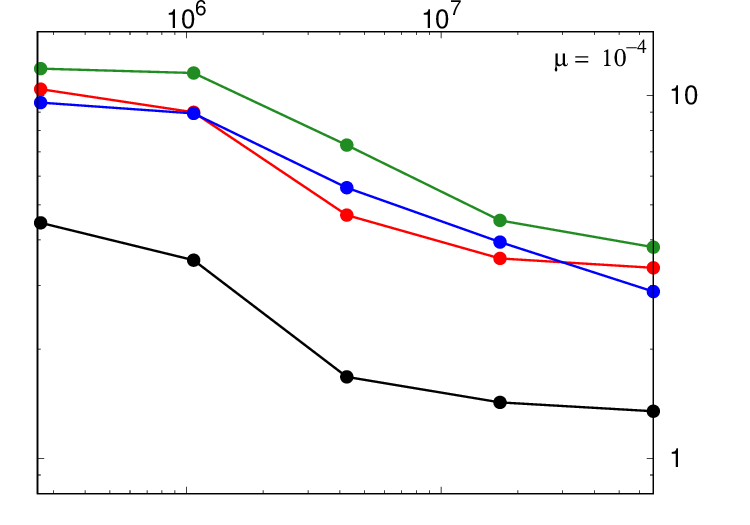}
  \includegraphics[height=0.24\textwidth,trim={0 0 44 10},clip]{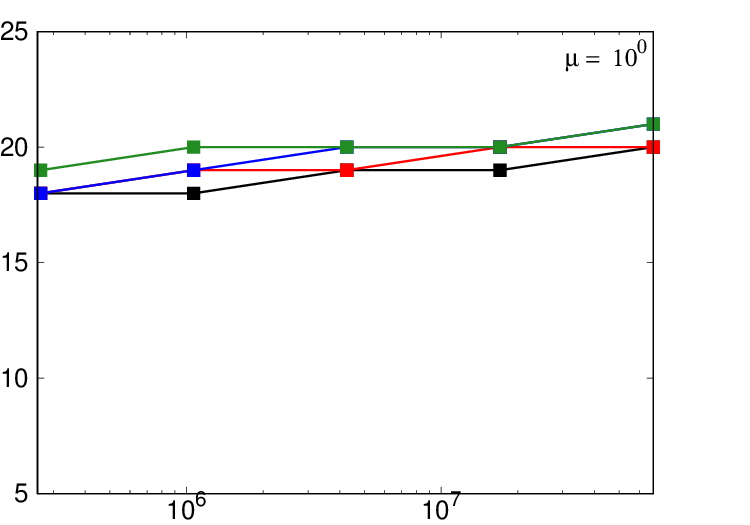}
  \includegraphics[height=0.24\textwidth,trim={16 0 44 10},clip]{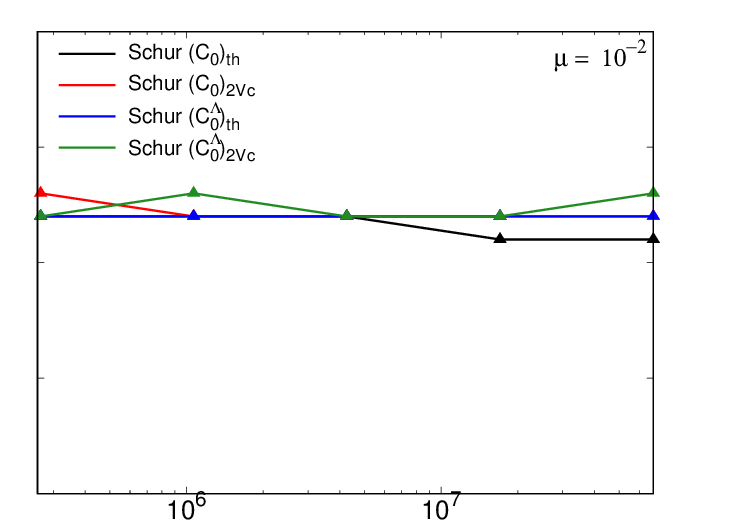}
  \includegraphics[height=0.24\textwidth,trim={16 0 0 10},clip]{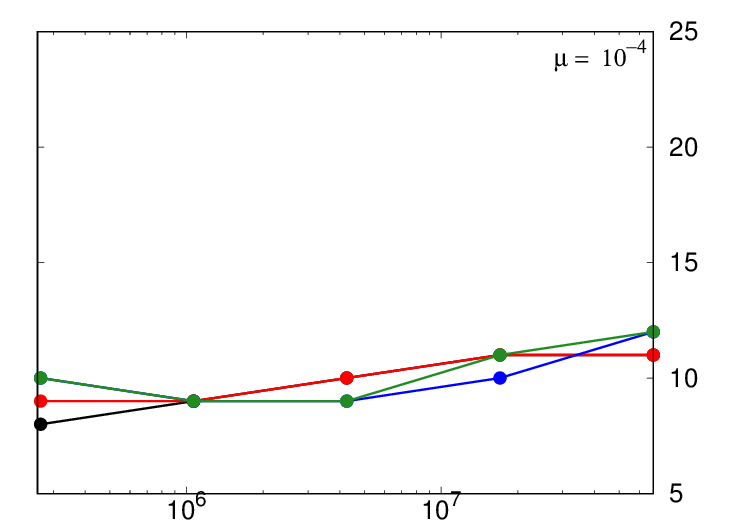}
  \caption{Comparison of the four  preconditioners
    $(C_0)\lo{th}$, $(C_0)\lo{2Vc}$,
    $(C_0^\Lambda)\lo{th}$, and $(C_0^\Lambda)\lo{2Vc}$.  Y-axis:
    throughput in \bal{kdofs per s} for the top row and number of GMRES
    iterations for the bottom row. X-axis: total number of dofs. Left column: $\mu=1.$ Center column: $\mu=10^{-2}.$
     Right column: $\mu=10^{-4}.$}\label{fig:efficiency_schur_l0}%
  \vspace{-\baselineskip}
\end{figure}

We report in the \bal{top panels} of Figure~\ref{fig:efficiency_schur_l0}
the throughput (kdofs per s) as a function of the total number of degrees of freedom
for the four preconditioners in \eqref{Cahouet_Chabard_precondiotioners} and for three values of the
viscosity: $\mu\in\{1,10^{-2},10^{-4}\}$.  We observe that the three
preconditioners $(C_0)\lo{2Vc}$, $(C_0^\Lambda)\lo{th}$, and
$(C_0^\Lambda)\lo{2Vc}$ behave similarly. The preconditioner
$(C_0)\lo{th}$ is slower because the mass matrix problem
$M_\bV\sfW=\sfZ$ is solved almost exactly (up to $10^{-10}$ accuracy),
which is a waste of resource since these operations are only used for
preconditioning purposes.

We clearly observe that the four methods can be ordered in throughput performance.
The best method is $(C_0^\Lambda)\lo{2Vc}$, then comes  $(C_0)\lo{2Vc}$ followed by  $(C_0^\Lambda)\lo{th}$,
and $(C_0)\lo{th}$ is the slowest method.

Depending on $\mu$, the throughput for $(C_0^\Lambda)\lo{2Vc}$ is in the range [6,10]$\CROSS10^3$kdof/s on
the coarsest grid and in the range
[2,4]$\CROSS10^3$kdof/s on the finest grid, which is on par with (or noticeably better than) what is
reported in the literature (see
Table~\ref{tab:litterature_comparisons}).  The weak scaling is
acceptable as the throughput  only decreases by a factor 3 when the
number of dofs is multiplied by $255$, (recall that the range of the number of dofs
is [267,427, 68,076,547]).  This behavior is coherent with the slight
loss of scaling observed in the inversion of $B \Lambda_V^{-1}B$
documented in \S\ref{Sec:B_LambdaV_Bt} which is likely due to communications between
nodes.

We show in the \bal{bottom panels} of Figure~\ref{fig:efficiency_schur_l0} the
number of outer GMRES iteration as a function of the total number of
degrees of freedom for the four preconditioners in \eqref{Cahouet_Chabard_precondiotioners} and for the same
values of the viscosity: $\mu\in\{1,10^{-2},10^{-4}\}$.  We observe
that all the methods are very robust with respect to the mesh size and
the viscosity, including $(C_0)\lo{th}$. But, although all the methods are optimal in terms of
GMRES iterations, there is a slight loss of weak scalability on the
efficiency due to communications between nodes. \bal{The method that has the best throughput is $(C_0^\Lambda)\lo{2Vc}$.}

\bal{This series of tests shows again that just looking at the number of
outer GMRES iterations is not fully informative.  The throughput is
the most important differentiating factor between preconditioners.}

\begin{remark}[Replacing $BM_q^{-1}B\tr$ by $L_Q^{-1}$]
  We finish this section by saying a few words regarding the following
  alternative preconditioner also mentioned in
  \cite{Cahouet_Chabard_1988} \bal{which consists of replacing
  $BM_q^{-1}B\tr$ by $L_Q^{-1}$ in
  \eqref{Cahouet_Chabard_precondiotioners}:}
 \begin{align}
    (C_0^\Delta)\lo{th}& \eqq \mu  (M_Q)\lo{th}^{-1} + \dt^{-1} (L_Q)\lo{th}^{-1} .
\end{align}
We show in figure~\ref{fig:residuals_C2}, the GMRES residual reported
by PETSc as a
function of the number of GMRES iterations for the same selection of
meshes and viscosities as above.  We clearly observe that the
performance of this preconditioner is not robust with respect to the
number of degrees of freedom and the viscosity. We are not going to consider this
preconditioner in the rest of the paper \bal{and advise against using it}.

\begin{figure}[h]\centering
  \includegraphics[width=0.45\textwidth]{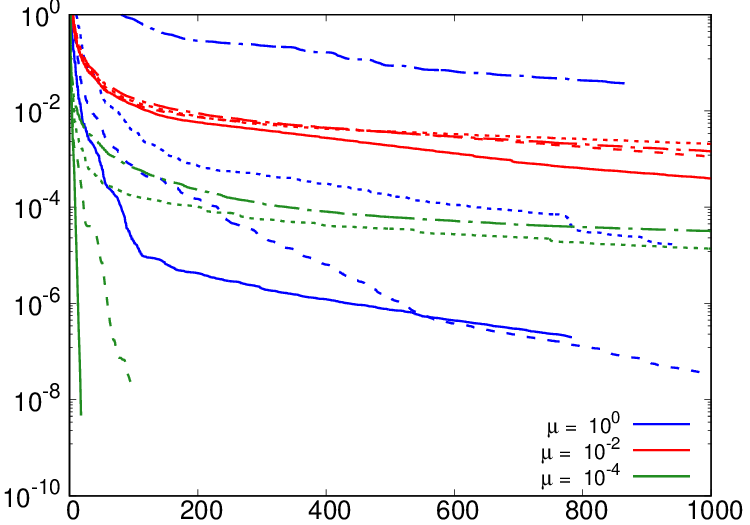}
  \caption{GMRES residual vs. iteration count for the preconditioner
    $(C_0^\Delta)\lo{th}$
      with $\mu=1$ ({\color{blue}
      blue}), $\mu=10^{-2}$ ({\color{red} red}), and $\mu=10^{-4}$
    ({\color{OliveGreen} green}). Number of velocity grid points: 118,785
    (\full); 473,857 (\dashed), 1,892,865 (\dotted) ; 7,566,337
    (\fulldot).}\label{fig:residuals_C2}%
  \vspace{-\baselineskip}
\end{figure}
\end{remark}

\subsection{Tests with augmented Lagrangian}

We now test the Schur complement technique with the augmented Lagrangian. We
solve \eqref{Schur_LA_problem} with the preconditioner
\eqref{AL_preconditioner} where $M_\bV$ is replaced by $\Lambda_\bV$, \ie
\begin{equation}
  (C_{\lambda}^\Lambda)\lo{2Vc}\eqq  \mu(1+\lambda) (M_Q)\lo{2Vc}^{-1} + \dt^{-1} (B\Lambda_\bV^{-1} B\tr)\lo{2Vc}^{-1}.
\end{equation}

We test the method with $\polP_2/\polP_1$ elements on the same meshes
as above and with $\lambda\in \{1,10\}$ and compare the results with
those obtained with the preconditioner $(C_{0}^\Lambda)\lo{2Vc}$
(which we recall corresponds to setting $\lambda=0$).

\begin{figure}[h]
  \centering
  \includegraphics[height=0.24\textwidth,trim={0 10 44 0},clip]{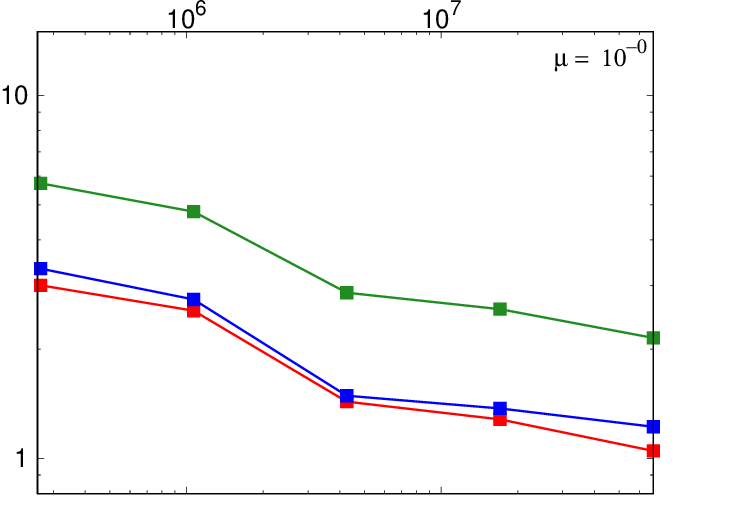}
  \includegraphics[height=0.24\textwidth,trim={16 10 44 0},clip]{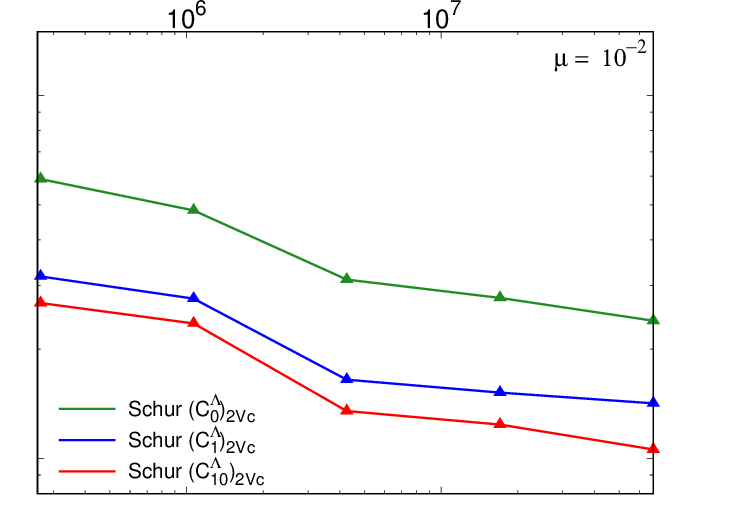}
  \includegraphics[height=0.24\textwidth,trim={16 10 0 0},clip]{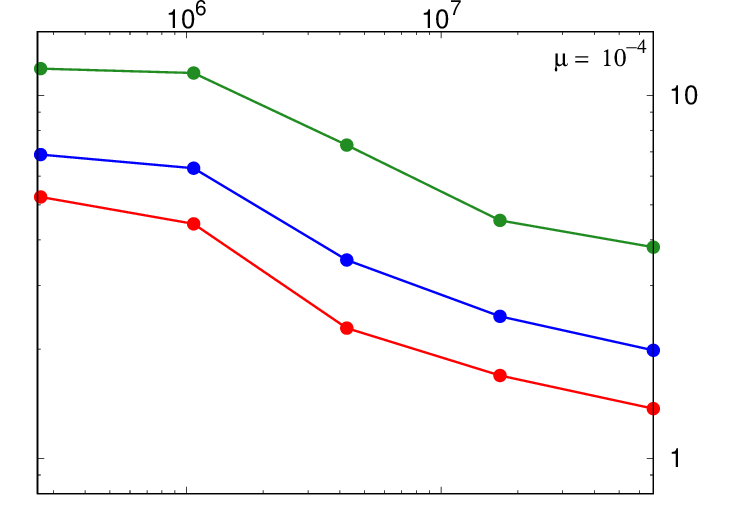}
  \includegraphics[height=0.24\textwidth,trim={0 0 44 10},clip]{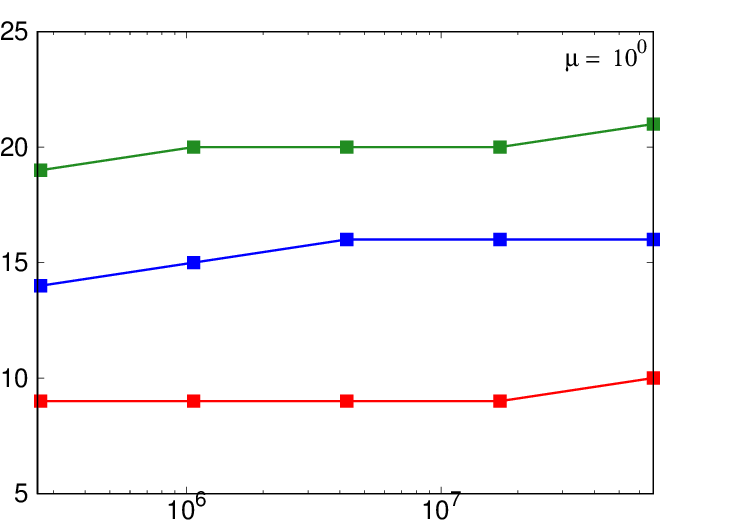}
  \includegraphics[height=0.24\textwidth,trim={16 0 44 10},clip]{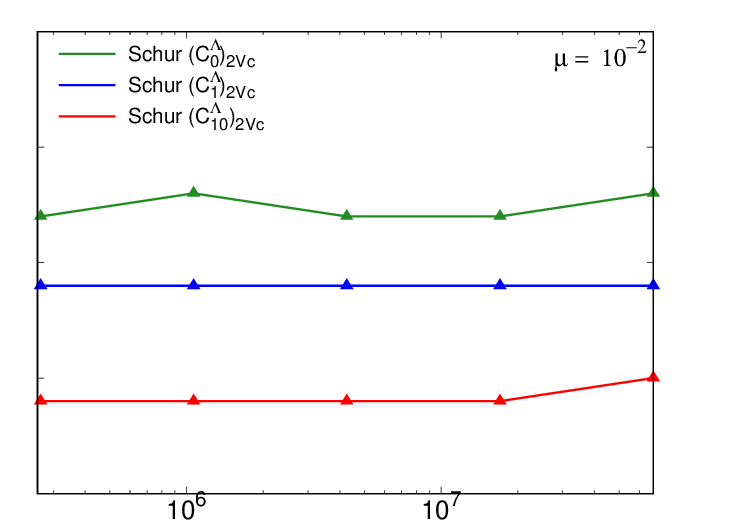}
  \includegraphics[height=0.24\textwidth,trim={16 0 0 10},clip]{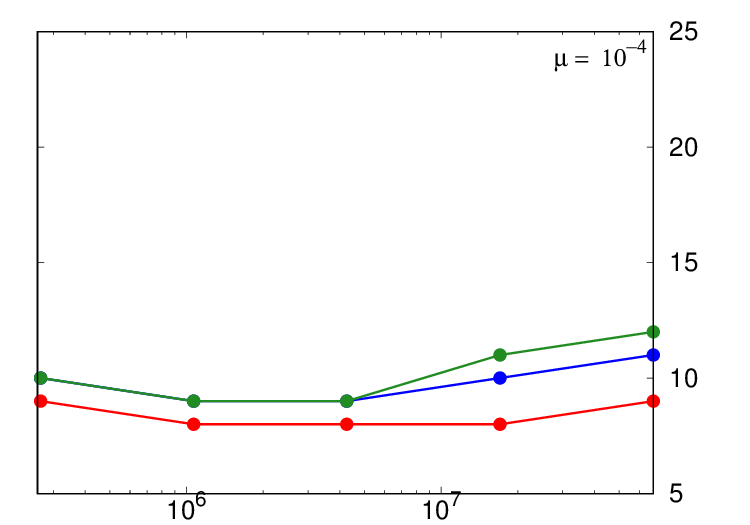}
  \caption{Comparison of  the three  preconditioners
    $(C_{0}^\Lambda)\lo{2Vc}$, $(C_{1}^\Lambda)\lo{2Vc}$,
    $(C_{10}^\Lambda)\lo{2Vc}$.  Y-axis:
    throughput in kdofs per s for the top row and number of GMRES
    iterations for the bottom row. X-axis: total number of dofs. Left column : $\mu=1.$
    Center column: $\mu=10^{-2}.$ Right column: $\mu=10^{-4}.$}\label{fig:efficiency_schur_lvarious}%
  \vspace{-\baselineskip}
\end{figure}

\bal{The throughput for the three methods $(C_{0}^\Lambda)\lo{2Vc}$, $(C_{1}^\Lambda)\lo{2Vc}$,
    $(C_{10}^\Lambda)\lo{2Vc}$ are shown in the top panels} of 
Figure~\ref{fig:efficiency_schur_lvarious}.  We observe that the three
methods $(C_{0}^\Lambda)\lo{2Vc}$, $(C_{1}^\Lambda)\lo{2Vc}$, and
$(C_{10}^\Lambda)\lo{2Vc}$ have similar weak scalability behavior. The throughput
decreases by a factor 3 when the total number of freedoms
is multiplied by 255 irrespective the value of $\lambda$. But quite surprisingly, we observe that the
throughput significantly decreases as $\lambda$ increases. 
\bal{The method that has the best throughput is $(C_{0}^\Lambda)\lo{2Vc}$.}

The  number of
outer GMRES iterations for the three methods is reported in the bottom panels of 
Figure~\ref{fig:efficiency_schur_lvarious}.
We observe that the number of GMRES iteration decreases as
$\lambda$ increases, as advertised in the augmented Lagrangian
literature.

\bal{The above series of tests show that the gain in GMRES iterations
obtained by increasing $\lambda$ is completely offset} by the
difficulty of solving the augmented Lagrangian velocity problem (as
reported numerous time in the literature). As a result, the time spent
solving the augmented Lagrangian velocity problem is not even compensated
by the reduction in GMRES iterations.

\subsection{Conclusions of \S\ref{Sec:Numerical_illustrations}}

We conclude that the preconditioner \bal{$(C_{0}^\Lambda)\lo{2Vc}$
  is} excellent. We also conclude that $(C_{0}^\Lambda)\lo{2Vc}$
outperforms in throughput
its augmented Lagrangian counterparts. Although the method
is robust in terms of GMRES iteration, a slight loss of weak
scalability is observed due to communications between nodes.

\section{Numerical illustration of the performance of Method~2} \label{Sec:Numerical_performance_Method2}

In this section we illustrate the performance of the solution method
described in \S\ref{Sec:Method_2}.

\subsection{Numerical details}

We recall that Method~2 consists of iteratively solving the coupled problem
\begin{equation}
 \polA_\lambda \begin{pmatrix} \sfU \\ \sfP\end{pmatrix}
 = \begin{pmatrix} \sfF \\ 0\end{pmatrix},\quad\text{with}\quad
 \polA_\lambda\eqq \begin{pmatrix} A_\lambda & -B\tr \\ B & 0 \end{pmatrix}.
\end{equation}
We use GMRES to solve this problem. For all the tests reported in the
paper, we keep iterating until the threshold  on the
relative residual reaches $10^{-10}$.  The following operations are
  done at each GMRES iterations: 
\begin{itemize}
    \item One matrix-vector multiplication with the matrix $\polA_\lambda$.
      \begin{itemize}
        \item If $\lambda=0$, all the operations are trivial.
        \item If $\lambda>0$, the augmented Lagrangian part of
          $A_\lambda$ requires solving one pressure mass problem.
          This is done using $(M_Q)\lo{th}^{-1}$.
    \end{itemize}
    \item One application of the preconditioner. This entails the following operations:
    \begin{itemize} 
    \item Two inversion of $\tA_{\lambda}$. This is done with CG with
      rel. thr. $10^{-10}$ using one of the preconditioners
      $(\tA_{\lambda,2})\lo{th}^{-1}$, $(\tA_{\lambda,2})\lo{2Vc}^{-1}$,
      $(\tA_{\lambda,3})\lo{th}^{-1}$, or $(\tA_{\lambda,3})\lo{2Vc}^{-1}$.
    \item Inversion of the Schur complement using either
      $(C_{\lambda}^\Lambda)\lo{th}$ or
      $(C_{\lambda}^\Lambda)\lo{2Vc}$ (see \S\ref{Sec:notation} and
      \S\ref{Sec:Numerical_details} for details).
    \end{itemize}
\end{itemize}
The method is tested using $\polP_2/\polP_1$ continuous Lagrange
elements on the same meshes as above and the viscosities
$\mu\in\{1,10^{-2},10^{4}\}$.

\subsection{Tests without augmented Lagrangian}

We start by investigating the performance of the method without the
augmented Lagrangian ($\lambda=0$).  We test the following four pairs:
\begin{subequations}
  \begin{align}
  ((C_{0}^\Lambda)\lo{th},(\tA_{0,2})\lo{th}^{-1}), \quad
  ((C_{0}^\Lambda)\lo{2Vc},(\tA_{0,2})\lo{2Vc}^{-1}), \\
  ((C_{0}^\Lambda)\lo{th},(\tA_{0,3})\lo{th}^{-1}), \quad
    ((C_{0}^\Lambda)\lo{2Vc},(\tA_{0,3})\lo{2Vc}^{-1}).
    \end{align}
\end{subequations}

\begin{figure}[h]
  \centering
  \includegraphics[height=0.24\textwidth,trim={0 10 44 0},clip]{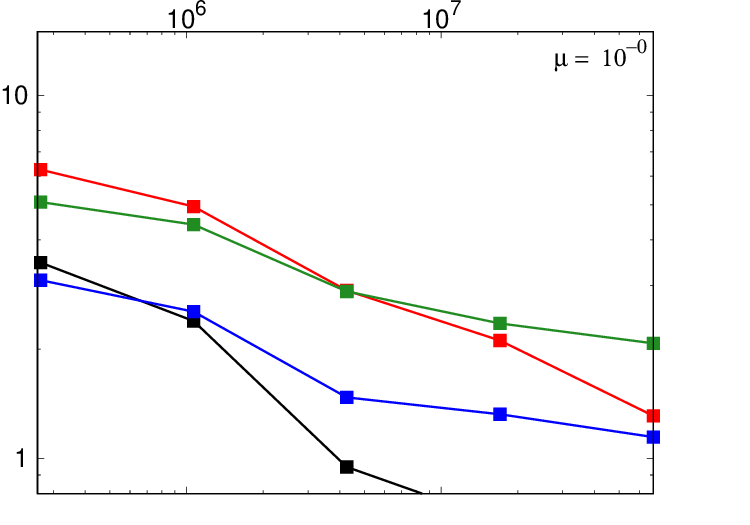}
  \includegraphics[height=0.24\textwidth,trim={16 10 44 0},clip]{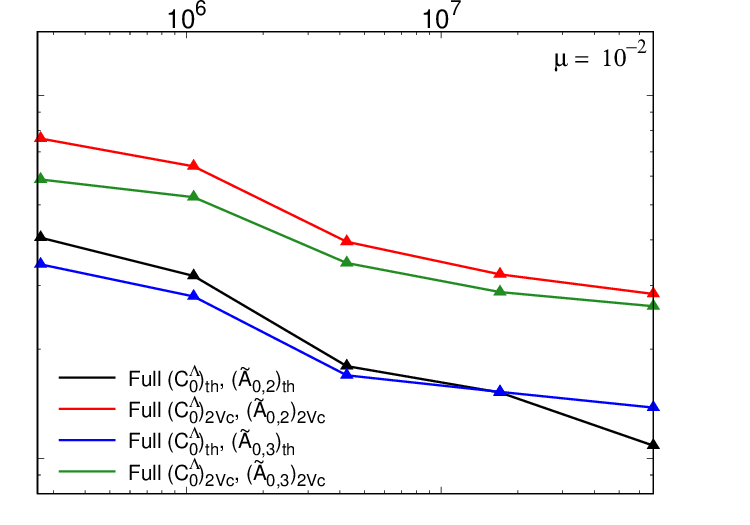}
  \includegraphics[height=0.24\textwidth,trim={16 10 0 0},clip]{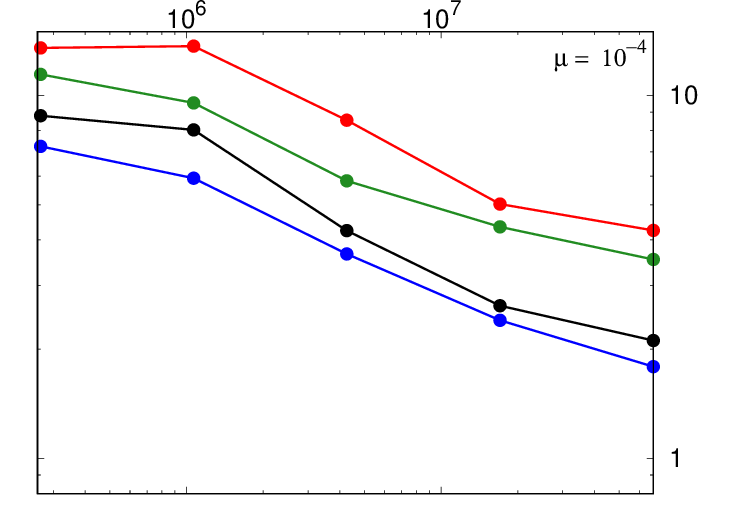}
  \includegraphics[height=0.24\textwidth,trim={0 0 44 10},clip]{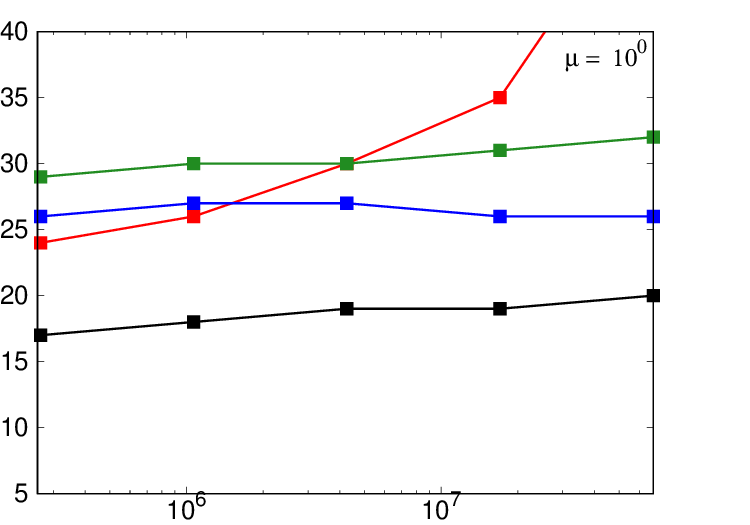}
  \includegraphics[height=0.24\textwidth,trim={16 0 44 10},clip]{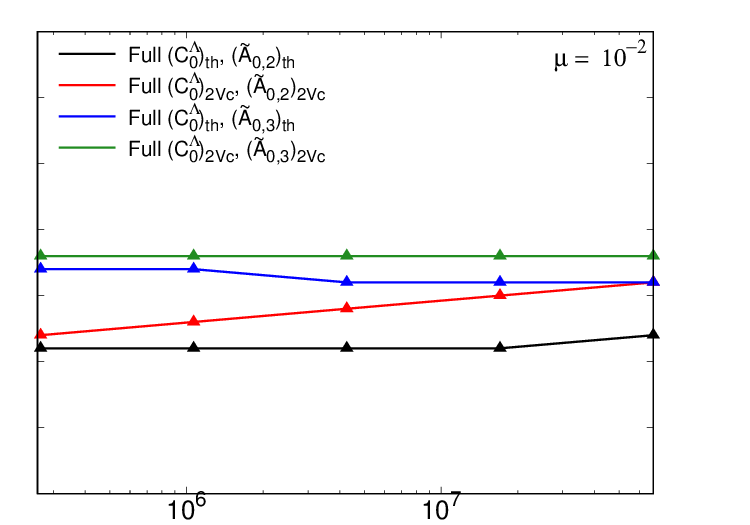}
  \includegraphics[height=0.24\textwidth,trim={16 0 0 10},clip]{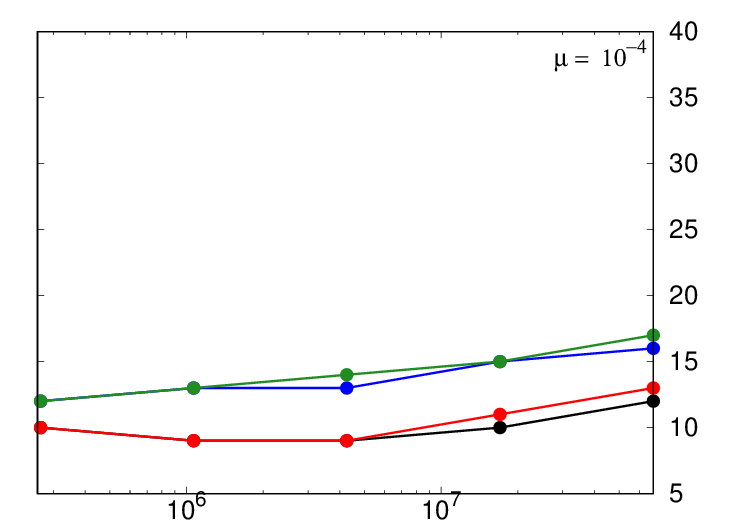}
  \caption{Comparison of the following four  preconditioner pairs for the full
    matrix $\polA_0$:
    $ ((C_{0}^\Lambda)\lo{th},(\tA_{0,2})\lo{th}^{-1})$,
    $((C_{0}^\Lambda)\lo{2Vc},(\tA_{0,2})\lo{2Vc}^{-1})$,
    $((C_{0}^\Lambda)\lo{th},(\tA_{0,3})\lo{th}^{-1})$, and
    $ ((C_{0}^\Lambda)\lo{2Vc},(\tA_{0,3})\lo{2Vc}^{-1})$.  Y-axis:
    throughput in kdofs per s for the top row, and number of GMRES
    iterations for the bottom row. X-axis: total number of dofs. Left column: $\mu=1.$ Center column: $\mu=10^{-2}.$
    Right column: $\mu=10^{-4}.$}\label{fig:efficiency_full_l0}%
  \vspace{-\baselineskip}
\end{figure}


We show the throughput in the top panels of Figure~\ref{fig:efficiency_full_l0} 
and the number of outer GMRES iterations in the bottom panels for these four methods.
We observe that only the two pairs
$((C_{0}^\Lambda)\lo{th},(\tA_{0,3})\lo{th}^{-1})$,
$((C_{0}^\Lambda)\lo{2Vc},(\tA_{0,3})\lo{2Vc}^{-1})$ behave well over
the entire range of viscosities. The pair
$((C_{0}^\Lambda)\lo{2Vc},(\tA_{0,3})\lo{2Vc}^{-1})$ seems to be the
most robust of the two in terms throughput.

The other two methods
$((C_{0}^\Lambda)\lo{th},(\tA_{0,2})\lo{th}^{-1})$,
$((C_{0}^\Lambda)\lo{2Vc},(\tA_{0,2})\lo{2Vc}^{-1})$ loose scalability
for $\mu=1$.  The pair
$((C_{0}^\Lambda)\lo{th},(\tA_{0,2})\lo{th}^{-1})$ performs well in terms
of outer GMRES iterations, but its low number of outer GMRES iterations does not compensate
for the difficulty of solving $(\tA_{0,2})\lo{th}^{-1}$. The other
method $((C_{0}^\Lambda)\lo{2Vc},(\tA_{0,2})\lo{2Vc}^{-1})$ has the
opposite problem \bal{when $\mu=1$.} Its number of outer GMRES iteration diverges when the
total number of degrees of freedom increases.

\begin{figure}[ht]\centering
    \includegraphics[width=0.45\textwidth]{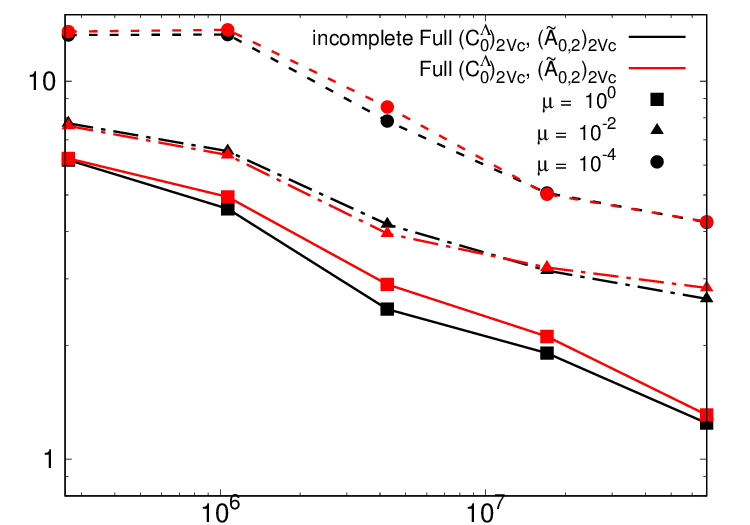}
    \includegraphics[width=0.45\textwidth]{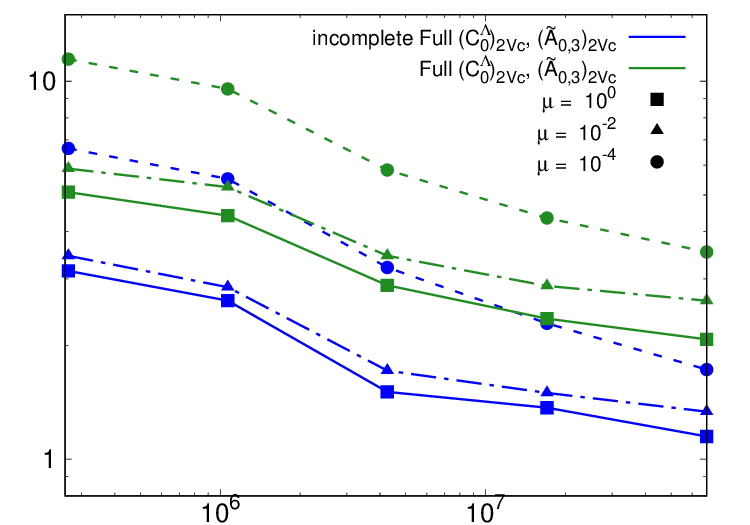}
  \caption{Left: ``incomplete Full''
  $(C^\Lambda_0)\lo{2Vc}, (\tA_{0,2})\lo{2Vc}$ and ``Full''
  $(C^\Lambda_0)\lo{2Vc}, (\tA_{0,2})\lo{2Vc}$. Right: ``incomplete Full''
  $(C^\Lambda_0)\lo{2Vc}, (\tA_{0,3})\lo{2Vc}$ and ``Full''
  $(C^\Lambda_0)\lo{2Vc}, (\tA_{0,3})\lo{2Vc}$. Y-axis:
    throughput in kdofs per s. X-axis: total number of dofs.}\label{fig:alternative_preconditioner}%
  \vspace{-\baselineskip}
\end{figure}
\begin{remark}[Alternative
  preconditioner] \label{Rem:alternative_preconditioner} \bal{It is
  sometimes advocated in the literature that a good preconditioner of
  $\polA_\lambda$ can be obtained by dropping the leftmost matrix on
  the right-hand side of \eqref{exact_factorization} as it requires one less inversion of $\tA_\lambda$, \ie
\begin{align*}
\polA_\lambda^{-1}\sim
\begin{pmatrix} \tA_\lambda^{-1}  & 0 \\ 0 & \tS_\lambda^{-1}\end{pmatrix}
\begin{pmatrix}I_\bV & 0 \\ -B \tA_\lambda^{-1} & I_Q\end{pmatrix}.
\end{align*}
We tested this approach using the pairs
$((C^\Lambda_0)\lo{2Vc}, (\tA_{0,2})\lo{2Vc})$ and
$((C^\Lambda_0)\lo{2Vc}, (\tA_{0,3})\lo{2Vc})$. We call these methods
``incomplete Full'' $(C^\Lambda_0)\lo{2Vc}, (\tA_{0,2})\lo{2Vc}$ and
``incomplete Full'' $(C^\Lambda_0)\lo{2Vc}, (\tA_{0,3})\lo{2Vc}$.  We
report the throughput for ``incomplete Full''
$(C^\Lambda_0)\lo{2Vc}, (\tA_{0,2})\lo{2Vc}$ in left panel of
Figure~\ref{fig:alternative_preconditioner}.  The throughput for
``incomplete Full'' $(C^\Lambda_0)\lo{2Vc}, (\tA_{0,3})\lo{2Vc}$ is
shown in the right panel. We observe that the throughput of
``incomplete Full'' $(C^\Lambda_0)\lo{2Vc}, (\tA_{0,2})\lo{2Vc}$ is
almost exactly the same as that of ``Full''
$(C^\Lambda_0)\lo{2Vc}, (\tA_{0,2})\lo{2Vc}$.  We also observe that
the throughput of ``incomplete Full''
$(C^\Lambda_0)\lo{2Vc}, (\tA_{0,3})\lo{2Vc}$ is significantly lower
than that of ``Full'' $(C^\Lambda_0)\lo{2Vc}, (\tA_{0,3})\lo{2Vc}$.
In conclusion, the benefits of this approach are not clear.}
\end{remark}

\subsection{Tests with augmented Lagrangian}


We now test the ``full'' method
with the augmented Lagrangian using $\lambda=1$ and compare the results with
the method without the augmented Lagrangian. We only consider the
following two methods:
  \begin{align}
  ((C_{1}^\Lambda)\lo{2Vc},(\tA_{1,2})\lo{2Vc}^{-1}), \quad
    ((C_{1}^\Lambda)\lo{2Vc},(\tA_{1,3})\lo{2Vc}^{-1}),
    \end{align}
    which we compare with the pairs  $((C_{0}^\Lambda)\lo{2Vc},(\tA_{0,2})\lo{2Vc}^{-1})$, 
    $((C_{0}^\Lambda)\lo{2Vc},(\tA_{0,3})\lo{2Vc}^{-1})$.
\begin{figure}[h]
  \centering
  \includegraphics[height=0.24\textwidth,trim={0 10 44 0},clip]{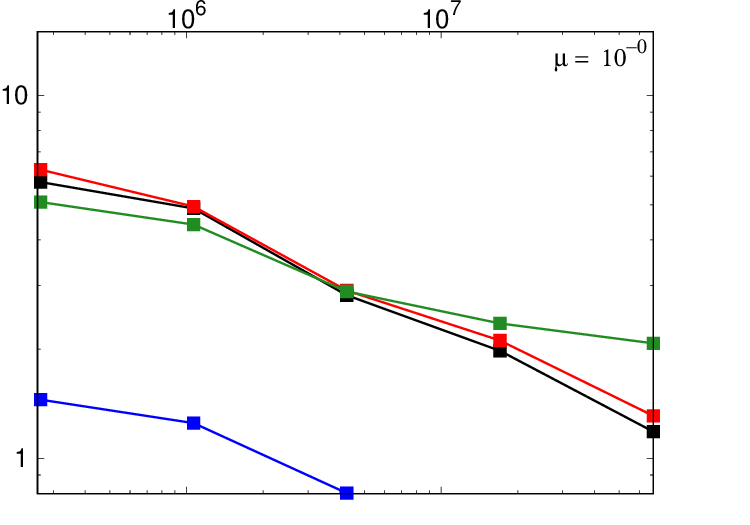}
  \includegraphics[height=0.24\textwidth,trim={16 10 44 0},clip]{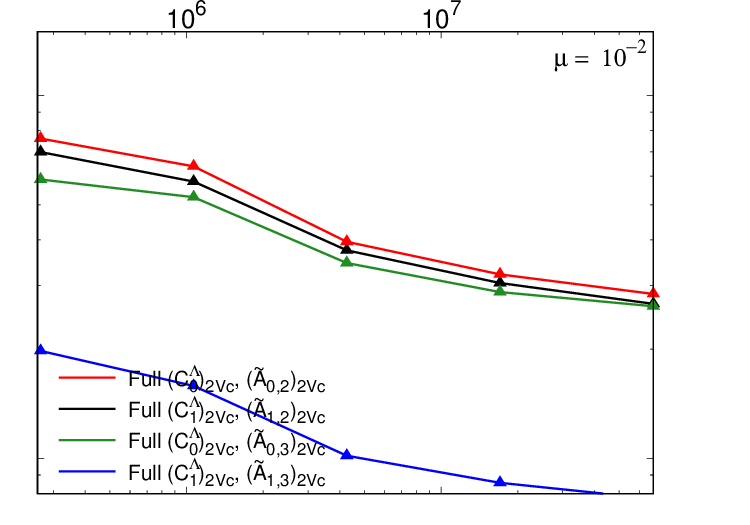}
  \includegraphics[height=0.24\textwidth,trim={16 10 0 0},clip]{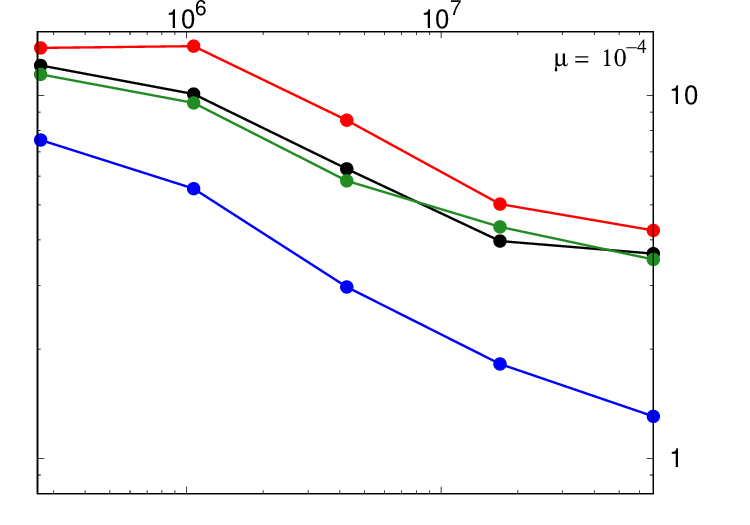}
  \includegraphics[height=0.24\textwidth,trim={0 0 44 10},clip]{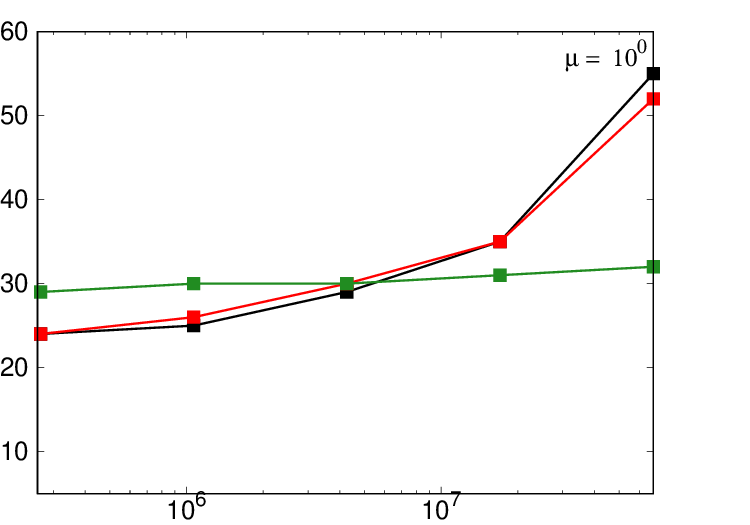}
  \includegraphics[height=0.24\textwidth,trim={16 0 44 10},clip]{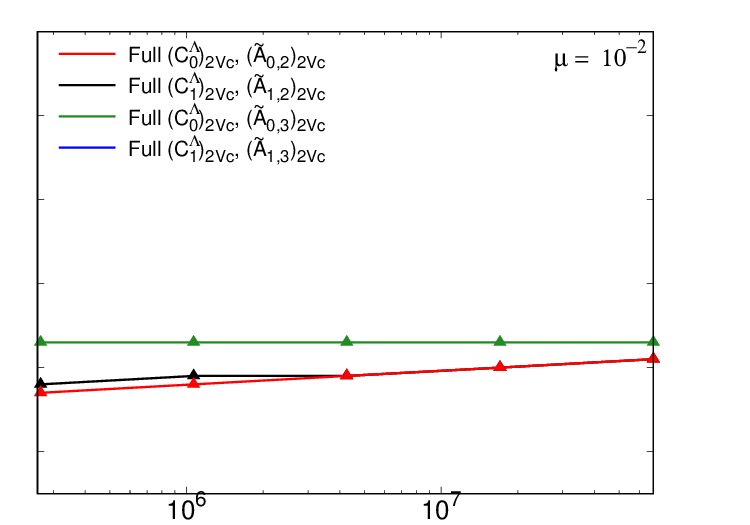}
  \includegraphics[height=0.24\textwidth,trim={16 0 0 10},clip]{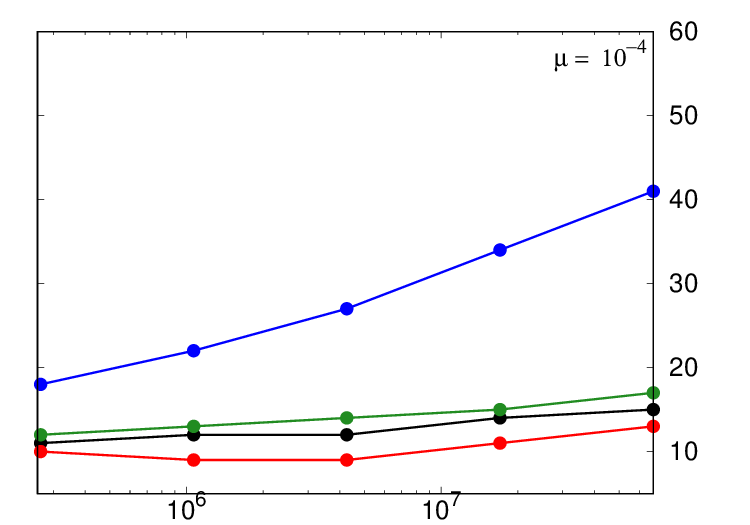}
  \caption{Comparison of
    $ ((C_{0}^\Lambda)\lo{2Vc},(\tA_{0,2})\lo{2Vc}^{-1})$,
     $((C_{1}^\Lambda)\lo{2Vc},(\tA_{1,2})\lo{2Vc}^{-1})$,
     $((C_{0}^\Lambda)\lo{2Vc},(\tA_{0,3})\lo{2Vc}^{-1})$, and
     $ ((C_{1}^\Lambda)\lo{2Vc},(\tA_{1,3})\lo{2Vc}^{-1})$.  Y-axis:
    throughput in kdofs per s for the top row and number of GMRES
    iterations for the bottom row. X-axis: total number of dofs. Left column: $\mu=1.$ Center column: $\mu=10^{-2}.$
    Right column: $\mu=10^{-4}.$}\label{fig:efficiency_full_lvarious}%
  \vspace{-1\baselineskip}
\end{figure}
    
The results of the simulations are shown in
Figure~\ref{fig:efficiency_full_lvarious}. The throughput
is reported in the top panels of the figure and the number
of GMRES iterations to achieve convergence is shown in the bottom
panels. We observe that the method
$((C_{1}^\Lambda)\lo{2Vc},(\tA_{1,3})\lo{2Vc}^{-1})$ is highly
inefficient ({\color{blue} blue} lines in the left and right panels of
Figure~\ref{fig:efficiency_full_lvarious})).  The other method
$((C_{1}^\Lambda)\lo{2Vc},(\tA_{1,2})\lo{2Vc}^{-1})$ is not robust
when $\mu=1$, which is somewhat surprising.  Here again we observe
that the difficulty of efficiently solving the augmented Lagrangian
velocity problem makes the methods not competitive (at least with the
tools we have used in the PETSc and BoomerAMG toolbox).

\section{Conclusions} \label{Sec:Conclusions}

We summarize the findings of the paper in this section. All the conclusions
recorded here hold assuming that all the subproblems are solved
iteratively (which is a necessity for very large
problems).  No direct solver is used at any stage. We also recall
that, to be representative of situations corresponding to the approximation
of the time-dependent Navier-Stokes equations,
all the tests done in the paper assume that the time step scales like
the mesh size, \ie $\dt =N^{-\frac12}$ in dimension 2, where $N$ is
the total number of velocity grid points (see \eqref{def_of_dt}).

\subsection{Cahouet\&Chabard}\label{Conclusion:Cahouet-Chabard}

We have shown that $(C_{0}^\Lambda)\lo{2Vc}$ and
$(C_{0}^\Lambda)\lo{th}$ are excellent preconditioners of the pressure
Schur complement matrix (see Figure~\ref{fig:efficiency_schur_l0}).

We have showed/confirmed that using the discrete pressure Laplacian
instead of the full matrix $BM_\bV^{-1} B\tr$ in the Cahouet\&Chabard
preconditioning of the pressure Schur complement is not a
robust strategy. This gives a methods that converges slowly on fine
meshes and when the thresholds is stringent (see
Figure~\ref{fig:residuals_C2}).  Hence, we recommend using
$BM_\bV^{-1} B\tr$ in the preconditioner even if this means solving
another linear problem using another iterative method.

\subsection{Augmented Lagrangian versus Cahouet\&Chabard}\label{Conclusion:AL_vs_CaCh}

We have confirmed that the preconditioned augmented
Lagrangian method requires less GMRES iterations than the traditional
Cahouet\&Chabard method for solving the pressure Schur complement
problem.  However, too much time is spent solving the velocity problem
in the augmented Lagrangian method.  The main difficulty consists of
inverting the discrete version of the grad-div operator.  Although the
augmented Lagrangian method requires less outer iterations to
reach convergence, the overall throughput of the augmented
Lagrangian method is significantly lower than that of the Cahouet\&Chabard method and
deteriorates as $\lambda$ grows (see
Figure~\ref{fig:efficiency_schur_lvarious}).

\subsection{Schur complement versus the full system} \label{Conclusion:Schur_vs_Full}


We have shown that, when the time step $\dt$
scales like the mesh size, there is no significant throughput gain 
by solving the full system instead of just solving the
pressure Schur complement.

\begin{figure}[h]
  \centering
  \includegraphics[height=0.24\textwidth,trim={0 10 44 0},clip]{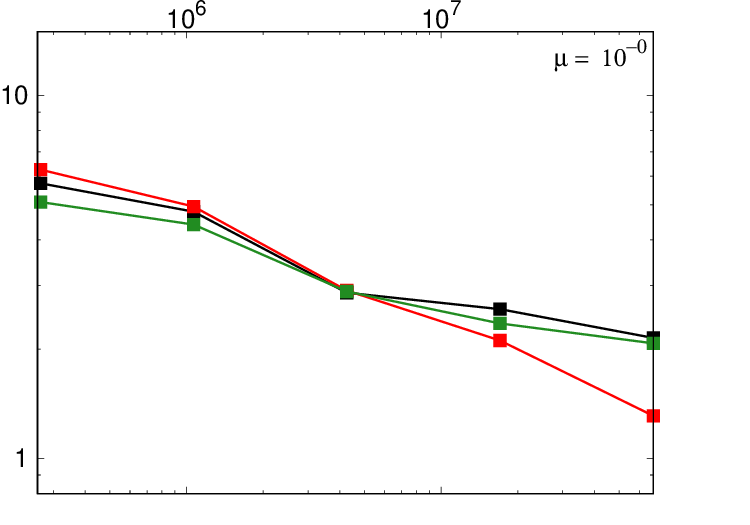}
  \includegraphics[height=0.24\textwidth,trim={16 10 44 0},clip]{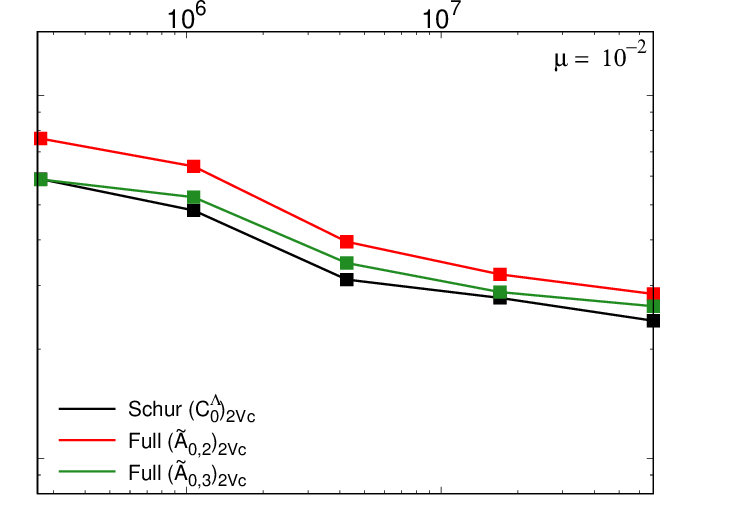}
  \includegraphics[height=0.24\textwidth,trim={16 10 0 0},clip]{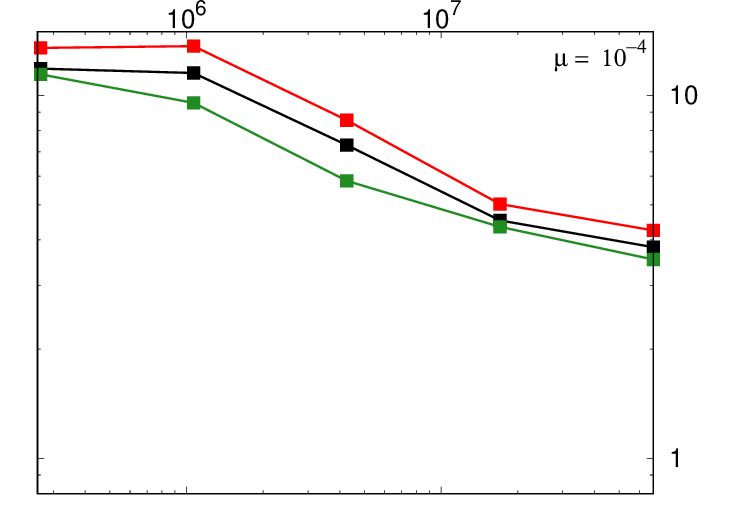}
  \caption{Comparison of the preconditioner $(C_{0}^\Lambda)\lo{2Vc}$
    for the pressure Schur complement with the preconditioner pairs
    ($(C_{0}^\Lambda)\lo{2Vc}, (\tA_{0,3}^\Lambda)\lo{2Vc}$) and
    ($(C_{0}^\Lambda)\lo{2Vc}, (\tA_{0,2}^\Lambda)\lo{2Vc}$) for the
    full system.  Y-axis: throughput in kdofs per s. X-axis: total
    number of dofs. Left panel: $\mu=1.$ Center panel: $\mu=10^{-2}.$
    Right panel: $\mu=10^{-4}.$}\label{fig:efficiency_best}%
  \vspace{-\baselineskip}
\end{figure}

We compare the two methods in Figure~\ref{fig:efficiency_best}. We
observe that when $\mu=1$, the Schur complement method has the best
throughput. For $\mu=10^{-2}$ and $\mu=10^{-4}$ the throughput
of the pressure Schur complement is marginally better than that of the
pair ($(C_{0}^\Lambda)\lo{2Vc}, (\tA_{0,3}^\Lambda)\lo{2Vc}$), and it
is marginally below  that of the pair
($(C_{0}^\Lambda)\lo{2Vc}, (\tA_{0,2}^\Lambda)\lo{2Vc}$). Overall the
two sets of methods are equivalent in terms of throughput.

\subsection{Comparisons with the literature}

We now compare the throughput that has been obtained in the paper
(solving either the pressure Schur complement or the full system) with
some results reported in the literature: \cite{Benzi_Wang_SISC_2011},
\cite{Farrell_Mitchell_Wechsung_SISC_2019},
\cite{Moulin_Jolivet_Marquet_CMAME_2019}, \cite{Larin2008ACS},
\cite{Shih_Stadler_Wechsung_SISC_2023}, \cite{Farell2021},
\cite{Niet_Wubs_IJNMF_2007}.

The comparisons are compiled in
Table~\ref{tab:litterature_comparisons}.
The column ``Source'' gives the reference and the table or figure we
refer to in this reference.  The column ``Mesh'' says which
discretization is used by the authors and also gives the space
dimension. The column ``Thr.'' gives the threshold on the relative
residual that is used in the reference. The column ``Nb. DoF'' gives the
range of total number of degrees of freedom explored in the reference.
The column ``Nb. Proc.'' gives the range of the number of processors
used for the tests reported in the reference. The column \bal{``TPS'' gives
the range of the throughput} achieved in the reference. The throughput
is rescaled to achieve $10^{-10}$ error on the relative
residual. For instance, if the number $\text{TPS}_{\text{ref}}$ is reported in
the reference with the relative tolerance $\text{thr}_{\text{ref}}$,
then we report in Table~\ref{tab:litterature_comparisons} the number
$\text{TPS}_{\text{ref}} \CROSS
\frac{|\log_{10}(\text{thr}_{\text{ref}})|}{10}$.

\begin{table}[h]\scriptsize\centering\addtolength{\tabcolsep}{-3.9pt}
  \begin{tabular}{|c|c|c|c|c|c|}\hline
  Source & Mesh & Thr. & Nb. DoF & Nb. Proc & TPS~(kdof/s) \\ \hline
  Fig. 9~\cite{Moulin_Jolivet_Marquet_CMAME_2019} &  3D $\polP_2/\polP_1$  &  $10^{-4}$  & 75M  & 256-2048 & $0.017$-$0.020$ \\  \hline
  Fig. 3~\cite{Shih_Stadler_Wechsung_SISC_2023} &  3D $[\polP_2 {\oplus}  B^F_3 ]^3/\polP_0$  & $10^{-6}$  &3.1M-1.6B  & 56-28672 & $0.21$-$0.56$    \\  \hline
  Fig. 6 and text~\cite{Farell2021}  &  3D $ [\polP_5]^3\polP^{\textup{disc}}_4 $  & $10^{-8}$  & 3.49M-30M   & 512-960 & $0.25$-$0.59$    \\  \hline
  Fig. 5.5~\cite{Farrell_Mitchell_Wechsung_SISC_2019} &  3D $[\polP_1 {\oplus}  B^F_3 ]^3/\polP_0$  & $10^{-5}$  &2.1M-1.1B   & 48-24576 & $0.56$-$0.67$   \\  \hline
  Fig. 6~\cite{Shih_Stadler_Wechsung_SISC_2023}  &  3D $[\polQ_2 ]^3/\polP_1^{\textup{disc}}$  & $10^{-6}$  &2.4M-151M   & 24-1536 & $0.77$-$1.37$    \\  \hline
  Fig. 7~\cite{Moulin_Jolivet_Marquet_CMAME_2019} &  2D $\polP_2/\polP_1$  &  $10^{-3}$  & 37k  & 1 & $0.091$\\ \hline
  Table V, VI~\cite{Larin2008ACS} &  3D $\polP_2/\polP_1$  &  $10^{-10}$  & 94k-786k  & 1 & $6.3$-$4.0$ \\  \hline
  Table 3~\cite{Chen_Jiao_ACM_2022} & 3D  $\polP_2/\polP_1$ & $10^{-6}$ & 263k-3.74M & 1 & $3.7$-$6.7$ \\ \hline
  Table II, VI, VII~\cite{Niet_Wubs_IJNMF_2007}  &  2D   & $10^{-6}$  & 768-197k   & 1 & $1.9$-$7.7$    \\  \hline
  Table 5.5~\cite{Farrell_Mitchell_Wechsung_SISC_2019} & 2D mesh & $10^{-6}$  & 160k  & 1 & $3.5$-$7.7$   \\ \hline
  Table 3.12~\cite{Benzi_Wang_SISC_2011} & 2D $\polQ_1/\polQ_1$ & $10^{-6}$  & 768-197k  & 1 & $3.0$-$9.1$   \\  \hline
  Table 3~\cite{Olshanskii_Zhiliakov_NLAP_2022} & 2D $\polP_2/\polP_1$ & $10^{-8}$ & 78k-1.3M & 1 & $9.1$-$12$ \\ \hline
  \rowcolor{gray!20}Present paper, Figure~\ref{fig:efficiency_best} & 2D $\polP_2/\polP_1$,  $\polP_3/\polP_2$ &  $10^{-10}$  & 267k-68M & 2-512 & $2.1$-$12$ \\ \hline
  \end{tabular}
  \caption{Comparison with the literature. Results from the paper are shown in the last row of the table in gray.}
\label{tab:litterature_comparisons}
\end{table}

We observe that the throughput obtained by our implementations
of the preconditioners for both the pressure Schur complement and the
full system (ranging from $2.1$ to $12$kdof/s to reduce the relative residual by 10 orders
of magnitude) compare favorably with the literature, thereby giving
credence to the claims made in the paper and summarized in \S\ref{Conclusion:Cahouet-Chabard} to
\S\ref{Conclusion:Schur_vs_Full} and in \S\ref{Conclusion:Schur_vs_Projection}.

\subsection{Schur complement versus projection methods}\label{Conclusion:Schur_vs_Projection}

Finally we want to compare the throughput of the methods discussed in
the paper with that of pressure-correction methods \`a la
Chorin--Temam.  We show in Table~\ref{tab:pred_corr_eff} the
throughput of a typical pressure-correction method which usually
consists of solving one velocity problem, one pressure Laplacian and
one pressure mass matrix per time step. The tests are done with the
preconditioners $(\tA_{0,3})\lo{2Vc}^{-1}$ for the velocity problem,
and the other two problems are solved using BoomerAMG. The threshold
on the relative residual is $10^{-10}$ for the three problems
(velocity problem, pressure Laplacian, pressure mass matrix).  We use
the same $\polP_2/\polP_1$ meshes as in
\S\ref{Sec:Numerical_illustrations} and
\S\ref{Sec:Numerical_performance_Method2}.

We observe that the throughput ranges from $278$kdof/s on the
coarsest mesh to $52$kdof/s on the finest mesh. These numbers have to
be compared to $12$kdof/s  and $2.1$kdof/s  for the methods discussed in the
paper (here we use the best numbers with the best methods carefully
optimized). The ratio of CPU time efficiency ranges from $23$ to $25$ in the best case scenario. In
conclusion, the methods discussed in the paper are on average 25 times
slower than pressure-correction methods.

\begin{table}[h]\scriptsize\centering
  \begin{tabular}{|c|c|c|c|c|c|c|}\hline
    \multirow{2}{*}{Visc. $\mu$} & Nb. dofs. & 267,427 & 1,066,499 & 4,259,587 & 17,025,539 & 68,076,547 \\
     & Procs & 2 & 8 & 32 & 128 & 512 \\ \hline
     \multirow{2}{*}{$10^{-0}$} & Times (s) & 1.046 & 1.308 & 2.296 & 3.543 & 3.665 \\
     & TPS (kdof/s)                & 128 & 102 & 56 & 38 & 37  \\ \hline
     \multirow{2}{*}{$10^{-2}$} & Times (s) & 1.033 & 1.313 & 2.162 & 3.624 & 3.055 \\
     & TPS (kdof/s)                & 129 & 102 & 62 & 36 & 44  \\ \hline
    \multirow{2}{*}{$10^{-4}$} & Times (s) & 0.479 & 0.586 & 1.061 & 2.779 & 2.568 \\
     & TPS (kdof/s)                 & 278 & 227 & 127 & 48 & 52  \\ \hline
   \end{tabular}
  \caption{Throughput for pressure-correction methods (kdof/s).}\vspace*{-\baselineskip}
  \label{tab:pred_corr_eff}
\end{table}

At the time of this writing, gaining a speedup factor of 25 on
preconditioned pressure
Schur complement techniques seems problematic
without a genuine algorithmic breakthrough. The augmented Lagrangian
method does not seem to be one of those. Of course one can argue
that Schur complement techniques are more accurate in time than
projection methods as they do not induce any time splitting
error. Moreover, they can be used with larger time steps as they can
be implemented with unconditionally stable time stepping technique.
More research on this topic is warranted.


\bibliographystyle{abbrvnat}
\bibliography{ref}
\end{document}